\documentclass[A4paper,10pt]{article}
\usepackage{amsmath,amssymb,mathrsfs}
\date{}

\def\func#1{\mathop{\rm #1}}
\renewcommand{\mathcal}{\mathscr}
\def\cl#1{{\mathscr #1}}
\newcommand{\dlines}{\displaylines}
\let\pa=\partial
\def\tr{\mathop{\rm tr}}
\newcommand{\field}[1]{\mathbb{#1}}
\newcommand{\R}{\field{R}}
\newcommand{\finedim}{\par\hfill$\blacksquare$\hphantom{$\blacksquare$$\blacksquare$}\par\noindent\ignorespaces}
\renewcommand{\th}{\theta}

\newcommand{\bT}{\field{T}}

\newcommand{\bS}{\field{S}}

\newcommand{\C}{\field{C}}

\newcommand{\Var}{{\rm Var}}

\newcommand{\e}{{\rm e}}

\newcommand{\cS}{{\mathcal S}}

\def\tin#1{\par\noindent\hskip3em\llap{#1\enspace}\ignorespaces}
\def\cl#1{{\mathcal #1}}

\def\cS{{\field S}}
\def\E{{\rm E}}

\def\P{{\rm P}}

\def\bC{{\field C}}

\def\paref#1{(\ref{#1})}

\newtheorem{theorem}{Theorem}[section]
\newtheorem{rema}[theorem]{Remark}

\newtheorem{assum}[theorem]{Assumption}
\newtheorem{cor}[theorem]{Corollary}
\newtheorem{definition}[theorem]{Definition}
\newtheorem{example}[theorem]{Example}

\newtheorem{lemma}[theorem]{Lemma}
\newtheorem{prop}[theorem]{Proposition}
\newtheorem{proposition}[theorem]{Proposition}
\newtheorem{remark}[theorem]{Remark}

\numberwithin{equation}{section}
\newcommand{\proof}{{\it Proof}.\ \/}
\begin{document}
\title{\huge\sc Fourier coefficients of invariant
random fields on homogeneous spaces of compact groups}
\author{P.~Baldi, S.~Trapani\\
{\sl Dipartimento di Matematica, Universit\`a
di Roma {\it Tor Vergata}, Italy}
}
\maketitle

\begin{abstract}
\noindent Let $T$ be a random field invariant under the action of
a compact group $G$. In the line of previous work
we investigate properties of the Fourier coefficients as orthogonality and Gaussianity.
In particular we give conditions ensuring that independence
of the random Fourier coefficients implies Gaussianity.
As a consequence, in general, it is not possible to simulate a
non-Gaussian invariant random field through its Fourier expansion
using independent coefficients.
\end{abstract}

\smallskip

\smallskip

\noindent{\it Acknowledgements}

\noindent The authors wish to thank the Referee for many illuminating remarks.
\section{Introduction}\label{intro}%
Recently much interest has been attracted to the investigation of
properties of  random fields on the sphere $\bS^2$ that are
invariant (in distribution) with respect to the action of the
rotation group $SO(3)$, highlighting a certain number of interesting
features (see \cite{dogiocam}, \cite{BM06} e.g.). This interest is
motivated mainly by the modeling and the investigation of
cosmological data.

For instance in \cite{BM06} it was proved that assumptions of independence of the Fourier coefficients of the
development in spherical harmonics
\begin{equation}\label{dev-SH}
T=\sum_{\ell=1}^\infty \sum_{m=-\ell}^\ell a_{\ell m}Y_{\ell m}
\end{equation}
of the real random field $T$, in addition to invariance, imply
Gaussianity. More precisely it was proved that if the coefficients
$a_{\ell m}$, $\ell=1,2,\dots, 0\le m\le \ell$ are independent, then
the field is necessarily Gaussian (the other coefficients are
constrained by the condition $a_{\ell,-m}=(-1)^m\overline {a_{\ell
m}}$). This result implies, in particular, the relevant consequence
that a non Gaussian invariant random field on $\bS^2$ {\it cannot}
be simulated using independent coefficients.

It is then a natural question whether this property also holds for
invariant random fields on more general structures. A result in this
direction was obtained in \cite{MR2342708} where it was proved that
a similar statement holds in general for an invariant random field
on the homogeneous space of a compact group, provided the
development is made with respect to a suitable Fourier basis
satisfying a particular condition (see Assumption \ref{assum0}
below).

The object of this paper is to pursue the line of investigation started in \cite{MR2342708}, in
particular in the direction of determining in which situations
Assumption \ref{assum0} is true. Actually we do not know whether
this condition holds in general and its investigation is one of the
main objects of this paper. Until now it is known to hold in the case of $\bS^2$ for the basis given by the spherical harmonics
$(Y_{\ell m})_m$, if $\ell>1$ (see \cite{dogiocam} p.145).

We shall give a new condition, equivalent to Assumption \ref{assum0}, which will enable us to prove that it is
satisfied for every self-conjugated Fourier basis of
the sphere $\bS^2$ (and not just the spherical harmonics) and also
for an important class of self-conjugated bases of the sphere $\bS^3$. In particular it is non possible to simulate 
a non Gaussian invariant real random field on $\bS^2$ using independent coefficients with respect to any self-conjugated 
bases of the irreducible $G$-module of $L^2(\bS^2)$.

Besides this characterization of Gaussianity, we discuss other properties of the
Fourier coefficients of an invariant
random field as orthogonality (it is well known that the $a_{\ell m}$'s of the
development \paref{dev-SH} on $\bS^2$ for an invariant random field
with finite variance are known to be orthogonal, see
\cite{dogiocam}, p.140) and invariance of their distribution with
respect to rotations of the complex plane.
\bigskip

\noindent The plan of the paper is as follows. In \S\ref{par-fourier} we
recall the main properties of the Fourier development of a random
field on the homogeneous space $\cl X$ of a compact group and give
necessary and sufficient  conditions for its invariance in terms of its development. In \S\ref{par-properties} we
investigate properties of its coefficients as orthogonality, among
other things. It turns out that, unlike the case of $\bS^2$, they are not orthogonal in general and precisions are made concerning this phenomenon. In \S\ref{par-complex} we recall (from
\cite{MR2342708}) results giving the characterization of Gaussianity
 which is our main concern. As mentioned above, these results, in many cases of
interest, hold under the assumption that  the Fourier basis that is
chosen for the Fourier development enjoys a certain property
(Assumption \ref{assum0}) with respect to the action of the group.
This section also contains a converse result, giving conditions on
Gaussian coefficients in order to produce an invariant random field.

The remainder of the paper is devoted to the investigation of the
validity of Assumption \ref{assum0}. In \S\ref{par-main} we give a new
equivalent condition which
is the main tool for the investigation of the two main examples
($\bS^2$ and $\bS^3$) which are the objects of
\S\ref{par-sphere-related} and \S\ref{par-spheres}. We are also able
to prove that independence of the Fourier coefficients entails
Gaussianity in some situations of interest in which it is known that Assumption \ref{assum0} does not hold (Theorem \ref{accadi}), in particular covering the case of the basis of the spherical harmonics $(Y_{\ell m})_m$ for $\ell=1$ (which completes the proof of the result of \cite{BM06}).

Finally \S\ref{par-conclusions} points out some open questions.
\section{A.s. square integrable random fields and Fourier developments}\label{par-fourier}
Throughout this paper $\mathcal X=G/K$ denotes the homogeneous space of a compact group $G$, $K$ being a closed subgroup. We
denote $x\to gx$, $g\in G$ the  action of $G$ and $dg$ and $dx$ respectively the Haar measure
of $G$ and the $G$-invariant measure on $\mathcal X$. We assume that both $dx$ and
$dg$ are normalized and have total mass equal to $1$ and we
write $L^2(\mathcal X)$ for $L^2(\mathcal
X,dx)$ and $L^2(G)$ for $L^2(G,dg)$. The spaces $L^2$ are spaces of {\it complex valued} square integrable functions.

Let us denote by $\widehat G$ the set of equivalence classes of irreducible representations of $G$ and let, for every $\sigma\in\widehat G$, $H^\sigma$ a Hilbert unitary $G$-module
of class $\sigma$ fixed from now on. For every $f\in L^2(G)$ let
$$
\widehat f(\sigma):=\sqrt{\dim\sigma}\int_G f(g)\sigma(g^{-1})\, dg\ .
$$
$\widehat f(\sigma)$ is a linear endomorphism of $H^\sigma$. The Peter-Weyl theorem (see \cite{MR1143783} or \cite{B-D} e.g.) states that
\begin{equation}\label{PW-op}
f(g)=\sum_{\sigma\in \widehat G}\sqrt{\dim\sigma}\,\tr\bigl( \widehat f(\sigma)\sigma(g)\bigr)\ .
\end{equation}
It is immediate that this is a good definition not depending, in particular, on the $G$-invariant scalar product that is considered on the $G$-module $H^\sigma$.

As soon as an orthonormal basis $h_1,\dots,h_m$, $m=\dim \sigma$, of $H^\sigma$ is chosen, one can define the matrix coefficients
$$
D^\sigma_{ij}(g)=\langle \sigma(g)h_j,h_i\rangle
$$
and the corresponding matrix entries of $\widehat f(\sigma)$:
$$
\widehat f(\sigma)_{ij}=\sqrt{\dim\sigma}\int_Gf(g)D^\sigma_{ij}(g^{-1})\, dg=\sqrt{\dim\sigma}
\int_Gf(g)\overline{D^\sigma_{ji}(g)}\, dg=\sqrt{\dim\sigma}\,\langle f,D^{\sigma}_{ji}\rangle_{L^2(G)}
$$
so that \paref{PW-op} becomes
\begin{equation}\label{PW-coeff}
f(g)=\sum_{\sigma\in\widehat G}\sqrt{\dim\sigma}\sum_{i,j=1}^{\dim \sigma}{\widehat f(\sigma)_{ij}}\,D^{\sigma}_{ji}(g)=\sum_{\sigma\in\widehat G}{\dim\sigma}\sum_{i,j=1}^{\dim \sigma}\langle f,D^{\sigma}_{ji}\rangle_{L^2(G)}\,D^{\sigma}_{ji}(g)
\end{equation}
The Peter-Weyl theorem also entails that, for every orthonormal basis $h_1,\dots,h_m$, $m=\dim \sigma$ of $H^\sigma$, the normalized matrix elements $\sqrt{\dim \sigma}\,D^\sigma_{ij}$, $\sigma\in\widehat G$, $1\le i,j\le\dim \sigma$, form an orthonormal complete basis of $L^2(G)$, so that \paref{PW-coeff} is just the corresponding Fourier development and $f(\sigma)_{ij}$ is the coefficient corresponding to the element $\sqrt{\dim \sigma}\,D^\sigma_{ji}$ of this basis.

Let us denote $L_g$ the left action of $G$ on $L^2(G)$, that is
$L_gf(h)= f(g^{-1}h)$. It is immediate that the functions
$(D^\sigma(g)_{ij})_{1\le i\le \dim\sigma}$, appearing in the columns of $D^\sigma$, span a subspace of $L^2(G)$ that is invariant and irreducible with respect to this left action.

From \paref{PW-op} a similar development can be derived for $L^2(\cl
X)$, $\cl X=G/K$. Actually remark that, for a
fixed $x_0\in {\cl X}$, the relation $\widetilde f(g)=f(gx_0)$
uniquely identifies functions in $L^2({\cl X})$ as functions in
$L^2(G)$ that are right invariant under the action of the isotropy
group $K$ of $x_0$ (that is formed by the elements $k\in G$ such
that $kx_0=x_0$).
For such functions $f$ we have, for every $k\in K$,
$$
\widehat f(\sigma)=\!\int_G\!f(g)\sigma(g^{-1})\, dg=\!\int_G\!f(gk)\sigma(g^{-1})\, dg=\!\int_G\!f(t)\sigma(kt^{-1})\, dt=
\sigma(k)\!\int_G\!f(t)\sigma(t^{-1})\, dt=
\sigma(k)\widehat f(\sigma)
$$
which implies that $\widehat f(\sigma)$ is $H^\sigma_K$-valued, $H^\sigma_K$ denoting the subspace of $H^\sigma$ of vectors that are invariant under the action of $K$, i.e. $\widehat f(\sigma)\in End(H^\sigma,H^\sigma_K)$. Hence, for a choice of an orthonormal basis $h_1,\dots,h_m$ of $H^\sigma$ such that $h_1,\dots,h_k$ span $H^\sigma_K$, the matrix $\widehat f(\sigma)$ will have all zeros in the rows from the $(k+1)$-th to the $m$-th. For $f\in L^2(\cl X)$ we shall consider, for simplicity, that $\widehat f(\sigma)$ is a $\dim\sigma\times\dim\sigma$ matrix with zeros on every row but for the first $\dim H^\sigma_K$ ones, corresponding to the first elements of the basis, that are supposed to be $K$-invariant. Remark that $H^\sigma_K$ might be reduced to $\{0\}$.

As a consequence of the aforementioned Peter-Weyl theorem we have the decomposition, that we shall need later,
\begin{equation}\label{e-pw-group}
L^2(\cl X)=\bigoplus_{\sigma\in\widehat G} \bigoplus_{1\le i\le
\dim(H^\sigma_K)}V^\sigma_{i}
\end{equation}
where $\widehat G$ denotes the set of equivalence classes of
irreducible representations of $G$ and the $V^\sigma_{i}$ are irreducible $G$-modules of class $\sigma$.
\bigskip

\noindent We consider on $\mathcal X$ a real or complex random field
$(T(x))_{x\in \mathcal X}$. This means that there exists a
probability space $(\Omega,\cl F,\P)$ on which the r.v.'s $T(x)$ are
defined and we shall always assume joint measurability, i.e.
$(x,\omega)\to T(x,\omega)$ is $\cl B(\cl X)\otimes\cl F$
measurable, $\cl B(\cl X)$ denoting the Borel $\sigma$-field of $\cl X$.

$T$ is said to be {\it a.s. continuous} if the map $x\to T(x)$ is
continuous a.s. It is said to be {\it a.s. square integrable} if
\begin{equation}\label{eq-l2bound}
\int_{\mathcal X}|T(x)|^2\, dx<+\infty,\qquad a.s.
\end{equation}
Remark that a.s. square integrability does not imply existence of moments of the r.v.'s $T(x)$.
If $T$ is a a.s. square integrable random field on $G$ then the function $x\to T(x,\omega)$ belongs to $L^2(\cl X)$ a.s. and one can define ``pathwise''
\begin{equation}\label{stoch-PW}
\widehat T(\sigma)=\sqrt{\dim\sigma}\int_G T(g)\sigma(g^{-1})\, dg
\end{equation}
which is now a $End(H^\sigma)$-valued r.v. Similarly we have the analog of \paref{PW-op}, i.e.
\begin{equation}\label{PW-opstoch}
T(h)=\sum_{\sigma\in \widehat G}\sqrt{\dim\sigma}\,\tr\bigl( \widehat T(\sigma)\sigma(h)\bigr)
\end{equation}
or
\begin{equation}\label{PW-coeffstoch}
T(h)=\sum_{\sigma\in\widehat G}\sqrt{\dim\sigma}\sum_{i,j=1}^{\dim \sigma}{\widehat T(\sigma)_{ij}}\,D^{\sigma}_{ji}(h)=\sum_{\sigma\in\widehat G}{\dim\sigma}\sum_{i,j=1}^{\dim \sigma}\langle T,D^{\sigma}_{ji}\rangle_2\,D^{\sigma}_{ji}(h)
\end{equation}
the series converging a.s. in $L^2(G)$.

For a random field $T$ we define the rotated random field $T^g$ as $T^g(x)=T(gx)$.
\begin{definition}\label{def-iso1}
A a.s. square integrable random field $T$ on $\cl X$ is said to be
\emph{$G$-invariant} if, as a $L^2(\cl X)$-valued random variable,
it has the same distribution as the rotated random field $T^g$ for
every $g\in G$, in the sense that the joint laws of
\begin{equation}\label{invar-L2gen}
(T({x_1}),\dots,T({x_m}))\qquad\mbox{and}\qquad(T({gx_1}),\dots,T({gx_m}))
\end{equation}
coincide for every $g\in G$ and $x_1,\dots, x_m\in\cl X$.

More generally a family $(T_i)_{i\in\cl I}$ of random fields on $\cl X$ is said to be invariant if and only if for every choice of $g\in G$, $i_1,\dots,i_m\in\cl I$ and $x_1,\dots,x_m\in \cl X$, the joint laws of
\begin{equation}\label{invar-L2gen-multi}
(T_{i_1}(x_1),\dots,T_{i_m}(x_m))\quad\mbox{and}\quad (T_{i_1}(gx_1),\dots,T_{i_m}(gx_m))
\end{equation}
coincide.
\end{definition}
The following will have some importance later. We thank D.Marinucci and G.Peccati for informing us of the existence of this result.
\begin{prop}\label{fd-l2} Let $T$  a a.s. square-integrable invariant random field on $\cl X$ and define, for $f\in L^2(\cl X)$,
\begin{equation}\label{gen-coeff}
T(f):=\int_{\cl X} T(x) \overline{f(x)}\,dx
\end{equation}
Then, for every $g \in  G$
and every $f_1,\dots, f_m \in L^2(G)$, the two random variables
$$
(T(f_1),\dots,T(f_m))\quad \mbox{and}\quad
(T^g(f_1),\dots,T^g(f_m))
$$
have the same distribution.
\end{prop}
\proof (See \cite{mp-continuity}.
\begin{prop}\label{prop-fourier-invariance}
Let $T$ a a.s. square integrable random field on $G$. Then $T$ is invariant if and only if, for every $g\in G$, the two families of
r.v's
$$
(\widehat T(\sigma))_{\sigma\in\widehat G}\quad\mbox{and}\quad (\widehat T(\sigma)\sigma(g))_{\sigma\in\widehat G}
$$
are equi-distributed.
\end{prop}
\proof
Let us assume $T$ invariant and let $\sigma\in\widehat G$. Then for every $v,w\in H^\sigma$ the function $g\to\langle\sigma(g^{-1})v,w\rangle$
is bounded and therefore in $L^2(G)$. Therefore, according to Proposition \ref{fd-l2} and denoting by $\sim$ equality in law, we have for every $g\in G$,
$$
\dlines{
\langle T(\sigma)v,w\rangle=\sqrt{\dim\sigma}\int_G T(h)\langle \sigma(h^{-1})v,w\rangle\, dh\sim
\sqrt{\dim\sigma}\int_G T(gh)\langle \sigma(h^{-1})v,w\rangle\, dh=\cr
=\sqrt{\dim\sigma}\int_G T(t)\langle \sigma(t^{-1}g)v,w\rangle\, dt=
\sqrt{\dim\sigma}\int_G T(t)\langle \sigma(t^{-1})\sigma(g)v,w\rangle\, dt=\langle T(\sigma)\sigma(g)v,w\rangle\cr
}
$$
This being true for every $v,w\in H^\sigma$, we have that, as $End(H^\sigma)$-valued r.v.'s, $T(\sigma)$ and $T(\sigma)\sigma(g)$ have the same distribution. Quite similarly, only
in a just more complicated way to write,
$$
(T(\sigma_1),\dots,T(\sigma_n))\quad\mbox{and}\quad (T(\sigma_1)\sigma_1(g),\dots,T(\sigma_n)\sigma_n(g))
$$
have the same distribution as a $End(H^{\sigma_1})\oplus\dots \oplus End(H^{\sigma_n})$-valued r.v., thus proving the only if part of the statement. The converse follows easily from
development \paref{PW-opstoch}.
\finedim
Let $f\in L^2(\cl X)$ and $V\subset L^2({\cl X})$ an irreducible $G$-module. We can then consider the orthogonal
projection $P_Vf$ of $f$ on $V$. Similarly for
a a.s. square integrable random field $T$ on $\cl X$ let us denote $T_V$ its orthogonal projection on $V$.
Remark that by definition (the functions of $V$ are necessarily continuous)
$T_{V}$ is always a continuous random field.

Let us denote by $D_{ij}(g)$ the matrix elements of the left regular action of $G$ on $V$ with respect to the orthonormal basis $(v_i)_i$ of $V$ and let us consider the random coefficients of the development of $T$ with respect to this basis
\begin{equation}\label{def-coeff}
a_i=\int_{\cl X} T(x)\overline{v_i(x)}\, dx\ .
\end{equation}
We denote by $a$ the complex vector with components
$a_i$, $i=1,\dots, d$. Then the coefficients of the rotated random field $T^g$ are obtained through the relation
\begin{equation}\label{rotated=coefficients}
\begin{array}c
\openup 2pt\displaystyle
a^g_i=\int_{\cl X_{\phantom{g}}} T(gx)\overline{v_i(x)}\, dx=\int_{\cl X} T(x)\overline{v_i(g^{-1}x)}\, dx=\cr
\displaystyle =\sum_{k=1}^{d}\overline{D_{ki}(g)}\int_{\cl X} T(x)\overline{v_k(x)}\, dx=\sum_{k=1}^{d}\overline{D_{ki}(g)}T(v_k)=\sum_{k=1}^{d}{D_{ik}}(g^{-1})a_k
\end{array}
\end{equation}
that is
\begin{equation}\label{rotated2}
a^g=D(g^{-1})a
\end{equation}
As
$$
T_{V}(x)=\sum_{k=1}^{d} a_k v_k(x)\ ,
$$
it is immediate that $T_V$ is invariant if and only if the random vectors $a$ and $D(g)a$ have the same distribution.

With respect to the Peter Weyl decomposition \paref{e-pw-group} we have
$$
T=\sum_{\sigma\in\widehat G}\sum_{i=1}^{\dim H^\sigma_K}T_{V^\sigma_i}\ .
$$
Using the fact that the projectors are $G$-equivariant (i.e. commute with the action of $G$) it is easy to prove the following, not really unexpected, statement (anyway see Proposition 3 of \cite{Pepy} for a proof).
\begin{prop} $T$ is invariant if and only if the family $(T_{V^\sigma_i})_{\sigma\in\widehat G,1\le i\le \dim H^\sigma_K}$
of random fields is invariant.
\end{prop}
We shall therefore concentrate our attention mainly on the projected random fields $T_V$.
When dealing with a real random field it is natural to require that
the basis $v_1,\dots, v_d$ of the $G$-module $V$, with respect to which the coefficients are computed ``respects'' the real and imaginary parts and, in particular, if $V=\overline V$, that this basis is stable under conjugation.
As explained in \cite{MR2342708}, \S2 and the Appendix, it is actually possible to decompose $L^2(\cl X)$ into a direct sum of $G$-modules in the form
\begin{equation}\label{e.pw-dec2}
L^2(\cl X)=\bigoplus_{i\in {\cl I}^o}V_i\oplus
\bigoplus_{i\in {\cl I}^+}(V_i\oplus \overline {V_i})
\end{equation}
where the direct sums are orthogonal and
$$
i\in {\cl I}^o\Leftrightarrow V_i=\overline V_i,\qquad
i\in {\cl I}^+\Leftrightarrow V_i\perp\overline V_i\ .
$$
We can therefore choose an orthonormal basis $(v_{ik})_{ik}$ of $L^2(\cl X)$ such that ($d_i=\dim V_i$)

$\bullet$ for $i\in {\cl I}^o$, $(v_{ik})_{1\le k\le d_i}$ is
an orthonormal basis of $V_i$ stable under conjugation;

$\bullet$  for $i\in {\cl I}^+$, $(v_{ik})_{1\le k\le d_i}$ is
an orthonormal basis of $V_i$ and
$(\overline {v_{ik}})_{1\le k\le d_i}$ is an orthonormal basis
of $\overline V_i$.

It is immediate that if $T$ is a real random field and $i\in\cl I^o$ then $T_{V_i}$ is also a real random field. On the other hand, if $i\in {\cl I}^+$  then $T_{V_i}$ and $T_{\overline{V}_i}$ may not be real
(actually they cannot be real unless they vanish), whereas $T_{V_i}+T_{\overline{V_i}}$ will be real.
\begin{remark}\label{remarks2}\rm
Representations of a compact group $G$ are classically classified as of real, complex or quaternionic type (see \cite{B-D}, p. 93 e.g.). In order to be self-contained let us recall that a {\it conjugation} $J$ of a $G$-module $V$ is an antilinear ($J(\alpha v)=\overline \alpha J(v)$) equivariant map
$J:V\to V$.

A $G$-module $V$ is said to be {\it real} if there exists a conjugation $J:V\to V$ such that $J^2=1$ and
{\it quaternionic} if there exists a conjugation $J:V\to V$ such that $J^2=-1$. It is {\it complex} if it is neither real nor quaternionic.

The important thing is that an irreducible $G$-module is of one and only one of these types and that equivalent $G$-modules are necessarily of the same type. If an irreducible $G$-module $V\subset L^2(\cl X)$ is such that $\overline V=V$, the usual conjugation $J:v\to \overline v$ is a real conjugation, so that $V$ must be of real type.  In particular, if a representation is of quaternionic or complex type, it cannot contain in its isotypical space a $G$-module that is self-conjugated, so that in the decomposition \paref{e.pw-dec2} it cannot be of type $\cl I^o$.

The irreducible representations of even dimension of $SU(2)$ are quaternionic and the corresponding $G$-modules appearing in the Peter-Weyl decomposition of this group cannot, therefore, be self-conjugated.
\end{remark}
\section{Properties of the coefficients}\label{par-properties}
In this section we give results concerning two properties that are enjoyed by the coefficients $\widehat T(\sigma)_{ij}$, $\sigma\in\widehat G$, $1\le i,j\le\dim\sigma$, of the Fourier development of an invariant random field on $\cl X$.

A random field $T$ is said to have finite variance if
\begin{equation}\label{def-finite variance}
\E\Bigl(\int_{\cl X}|T(x)|^2\, dx\Bigr)<+\infty
\end{equation}
\begin{remark}\label{rem-finite-variance} \rm If \paref{def-finite variance} holds, then the map
$x\to T(x)$ necessarily belongs to $L^2(\cl X)$ a.s., so that $T$ is a.s. square integrable. Also, by the Cauchy-Scwartz inequality, if $T$ has finite variance, the random variables $T(f)$, $f\in L^2(\cl X)$, defined in \paref{gen-coeff}, have finite variance. In particular the Fourier coefficients of $T$, with respect to any Fourier basis, also have finite variance.
\end{remark}
It is well known (see \cite{BM06}, \cite{dogiocam} p.126) that in the case $\cl X=\mathbb S^2$, $G=SO(3)$, if $T$ is invariant and has finite variance, its Fourier coefficients with respect to the basis formed by the spherical harmonics are pairwise orthogonal. Our first concern in this section is to investigate this question in the case of a more general basis and for a general homogeneous space of a compact group.

Remark that, if $T$ is invariant and has finite variance, for every $\sigma$ that is not the trivial representation, for the matrix entries $\widehat T(\sigma)_{ij}$ we have
\begin{equation}\label{mean-zero}
\E[\widehat T(\sigma)]=\E[\widehat T(\sigma)\sigma(g)]=
\E[\widehat T(\sigma)]\int_G\sigma(g)\, dg=0
\end{equation}
\begin{theorem}\label{ortho-true} Let $T$ a finite variance invariant random field on $\cl X$ and $\sigma_1,\sigma_2\in\widehat G$.

a) If $\sigma_1$ and $\sigma_2$ are not equivalent, then, for every orthonormal bases of $H^{\sigma_1}$ and $H^{\sigma_2}$ the r.v.'s $\widehat T(\sigma_1)_{ij}$ and $\widehat T(\sigma_2)_{k\ell}$ are orthogonal, $1\le i,j\le \dim\sigma_1$, $1\le k,\ell\le \dim\sigma_2$.

b) If $\sigma_1=\sigma_2=\sigma$ let $\Gamma(\sigma)=\E[\widehat T(\sigma)\widehat T(\sigma)^*]$. Then $\mathop{\rm Cov}(\widehat T(\sigma)_{ij},\widehat T(\sigma)_{k\ell})=\delta_{j\ell}\Gamma(\sigma)_{ik}$. In particular coefficients belonging to different columns are orthogonal and the covariance between entries in different rows of a same column does not depend on the column.
\end{theorem}
\proof Recall first that $\widehat T(\sigma)_{ij}=\langle T,D^\sigma_{ji}\rangle_2$ so that thanks to Remark \ref{rem-finite-variance}  the rv.'s $\widehat T(\sigma)_{ij}$'s have themselves finite variance. Let us denote $D^\sigma_{ij}(g)$ the matrix elements of the action of $G$ on $H^\sigma$.

a)
As $(\widehat T(\sigma_1),\widehat T(\sigma_2))$ has the same joint distribution as $(\widehat T(\sigma_1)\sigma_1(g),\widehat T(\sigma_2)\sigma_2(g))$ for every $g\in G$ by Proposition \ref{prop-fourier-invariance},
$$
\dlines{
\E\big[\widehat T(\sigma_1)_{ij} \overline {\widehat T(\sigma_2)_{k\ell}}\big]=
\E\big[(\widehat T(\sigma_1)D^{\sigma_1}_{ij}(g)) \overline{( \widehat T(\sigma_2)D^{\sigma_2}_{k\ell}(g))}\,\big]=\cr
=\sum_{r=1}^{\dim\sigma_1}
\sum_{m=1}^{\dim\sigma_2}D^{\sigma_1}_{rj}(g)
\overline{D^{\sigma_2}_{m\ell}(g)}\E\big[\widehat T(\sigma_1)_{ir}\overline{\widehat T(\sigma_2)_{km}}\,\big]\ .\cr
}
$$
This being true for every $g\in G$, it is also true if we take the integral of the right hand-side over $G$ in $dg$. As the functions
$D^{\sigma_1}_{rj}$ and $D^{\sigma_2}_{m\ell}$ are orthogonal for every choice of the indices, the representations $\sigma_1$ and $\sigma_2$ being not equivalent, we find
$$
\E[\widehat T(\sigma_1)_{ij} \overline {\widehat T(\sigma_2)_{k\ell}}]=0\ .
$$

b) If $\sigma_1=\sigma_2=\sigma$, the previous computation gives
$$
\dlines{
\E\big[\widehat T(\sigma)_{ij} \overline {\widehat T(\sigma)_{k\ell}}\big]
=
\sum_{r,m=1}^{\dim\sigma}\E\big[\widehat T(\sigma)_{ir}\overline{\widehat T(\sigma)_{km}}\,\big]\int_GD^{\sigma}_{rj}(g)
\overline{D^{\sigma}_{m\ell}(g)}\, dg=\cr
=\sum_{r,m=1}^{\dim\sigma}\E\big[\widehat T(\sigma)_{ir}\overline{\widehat T(\sigma)_{km}}\,\big]\delta_{rm}\delta_{j\ell}=\delta_{j\ell}
\sum_{r=1}^{\dim\sigma}\E\big[\widehat T(\sigma)_{ir}\overline{\widehat T(\sigma)_{kr}}\,\big]=
\delta_{j\ell}\Gamma(\sigma)_{ik}\ .\cr
}
$$
\finedim
Theorem \ref{ortho-true} states that the entries of $\widehat T(\sigma)$ might not be pairwise orthogonal  and this happens when the matrix $\Gamma$ is not diagonal. This phenomenon is actually already been remarked by other authors (see \cite{bib:M} Theorem 2 e.g.).

Example \ref{exf-bijoux} below provides an instance of this phenomenon. Of course there are situations in which orthogonality is still guaranteed: when the dimension of $H^\sigma_K$ is one at most (i.e. in every irreducible $G$-module the dimension of the space $H^\sigma_K$ of the $K$-invariant vectors in one at most) as is the case for $G=SO(d)$, $K=SO(d-1)$, $G/K=\mathbb S^{d-1}$. In this case actually the matrix $\widehat T(\sigma)$ has just one row that does not vanish and $\Gamma(\sigma)$ is all zeros, but one entry in the diagonal.

In order to produce an example of invariant random field for which the entries of $\widehat T(\sigma)$ are not pairwise orthogonal, we recall first a well known definition.

Let $Z=Z_1+iZ_2$ a complex r.v. $Z$ is said to be {\it Gaussian complex valued} if $(Z_1,Z_2)$ is jointly Gaussian. $Z$ is said to be {\it complex Gaussian} if, in addition, $Z_1$ and $Z_2$ are independent and have the same variance. If $Z$ is centered this is equivalent to the requirement that their distribution is invariant with respect to rotations of the complex plane. We shall use the following properties.

$\bullet$ A centered Gaussian complex valued r.v. $Z$ is complex Gaussian if and only if $\E[Z^2]=0$.

$\bullet$ Two centered  complex valued Gaussian r.v.'s $Z_1$, $Z_2$ are independent if and only if $\E[Z_1\overline {Z_2}]=\E[Z_1{Z_2}]=0$.
\begin{example} \label{exf-bijoux} \rm Let $\sigma\in \widehat G$ and $V\subset L^2(G)$ an irreducible $G$-module of dimension $d\ge 2$ of class $\sigma$ and denote by the matrix $D^\sigma(g)$ the action of $G$ on $V$ with respect to a fixed basis. Let $Z_1,\dots,Z_d$ independent centered complex Gaussian r.v.'s such that $\E[|Z_j|^2]=1$ for every $j$. Let $B=(b_{ij})_{ij}$ the random matrix defined as $b_{ij}=\alpha_iZ_j$, $\alpha_i\in\bC$. Then the random field
$$
T(g)=\sqrt{d}\,\tr(BD^\sigma(g))
$$
is invariant and, as it is immediate that $\widehat T(\sigma)=B$, its coefficients $\widehat T(\sigma)_{ij}$ are not pairwise orthogonal. Let us check invariance. Let $C=\widehat T(\sigma)D^\sigma(g)$, then $c_{ij}=\alpha_i\sum_{k=1}^dZ_kD^\sigma_{kj}(g)=\alpha_iW_j$ where
$$
W_j=\sum_{k=1}^dZ_kD^\sigma_{k j}(g)
$$
In view of Proposition \ref{prop-fourier-invariance} we must therefore just prove that the $W_j$'s are complex Gaussian, independent
and that $\E[|W_j|^2]=1$. First it is immediate that they are Gaussian complex valued. We have also
\begin{equation}\label{exf1}
\begin{array}{c}
\displaystyle\E[W_j\overline {W_k}]=
\E\Bigl[\sum_{h,r=1}^dZ_h\overline{Z_r}
D^\sigma_{hj}(g)
\overline{D^\sigma_{rk}(g)}\Bigr]=
\sum_{h,r=1}^d\delta_{hr}D^\sigma_{hj}(g)
\overline{D^\sigma_{rk}(g)}=\\
=\displaystyle\sum_{r=1}^d
D^\sigma_{rj}(g)\overline{D^\sigma_{rk}(g)}=\sum_{r=1}^dD^\sigma_{k r}(g^{-1})D^\sigma_{rj}(g)=\delta_{kj}\ .
\end{array}
\end{equation}
Similarly, as $\E[Z_hZ_r]=0$ for every $1\le h,r\le d$ (recall that $\E[Z^2]=0$ for a centered complex Gaussian r.v. $Z$),
\begin{equation}\label{exf2}
\E[W_j{W_k}]=\E\Bigl[\sum_{h,r=1}^dZ_h{Z_r}
D^\sigma_{hj}(g)
{D^\sigma_{r\ell}(g)}\Bigr]=0\ .
\end{equation}
\paref{exf1} for $k=j$ gives $\E[|W_j|^2]=1$, whereas \paref{exf1} and \paref{exf2} together imply that $W_j$ and $W_k$, $k\not=j$, are independent. Finally \paref{exf2} for $k=j$ gives $\E[W_j^2]=0$ for every $j$ so that the $W_j$'s are complex Gaussian, which completes the proof.
\end{example}
\begin{cor}\label{explicative}
Let $T$ an invariant random field with finite variance on $\cl X$ and let $V\subset L^2(\cl X)$ an irreducible
$G$-module different from the constants. Then the coefficients $(a_k)_k$ of the development of the projection
$T_V$ of $T$ on $V$ with respect to any orthonormal basis of $V$ are centered, orthogonal, and have a common
variance $c$.
\end{cor}
\proof It is repetition of the arguments of the proof of b) of Theorem \ref{ortho-true}.
As pointed out in Remark \ref{rem-finite-variance} the coefficients $a_k$'s have themselves finite variance and, thanks to \paref{rotated=coefficients} and $V$ being different from the constants, they are also centered.
From \paref{rotated2} we have, for every $g\in G$,
$$
\E[a_k\overline{a_\ell}]=\E[(D(g)a)_k\overline{(D(g)a)_\ell}]=\sum_{j,r=1}^d D_{kr}(g)\overline{D_{\ell j}(g)}\E[a_r\overline{a_j}]\ .
$$
Integrating in $dg$ and using the orthonormality properties of the matrix elements $D_{ij}(g)$ we find
$$
\E[a_k\overline{a_\ell}]=\frac 1{\dim V}\sum_{j,r=1}^{\dim V}\delta_{k\ell}\delta_{rj}\E[a_r\overline{a_j}]=\frac 1{\dim V}\,\delta_{k\ell}\sum_{j=1}^{\dim V}\E[|a_j|^2]\ .
$$
For $k\not=\ell$ this gives immediately the orthogonality, whereas for $k=\ell$ we have
$$
\E[|a_k|^2]=\frac 1{\dim V}\sum_{j=1}^{\dim V}\E[|a_j|^2]
$$
so that the $a_k$'s have the same variance.
\finedim
Another feature appearing in the case $G=SO(3)$, $\cl X=\mathbb S^2$ is that the coefficients $a_{\ell m}$ of the development in spherical harmonics \paref{dev-SH} of an invariant random field have each a distribution that is invariant with respect to rotations of the complex plane {\it if} $m\not=0$. The following discussion aims to see what can be said in general concerning this property.
\begin{remark}\label{max-torus}\rm
Let $V\subset L^2(\cl X)$ be an irreducible $G$-module of dimension $d>1$ and $\bT\subset G$ a maximal torus. Let
\begin{equation}\label{max-decomposition}
V=\bigoplus_{k=1}^d U_k
\end{equation}
be a decomposition of  $V$ into orthogonal  irreducible components of the action of $\bT$ on $V$. As
$\bT$ is abelian, $\dim(U_k)=1$ for every $k=1,\dots,d$. Let $u_k\in U_k$ be a unit vector. Then
$L_tu_k=u_k(t^{-1}\cdot)=\chi_k(t) u_k$ for $t\in\bT$, where $\chi_k$ denotes the character of the representation of $\bT$ on $U_k$. If we consider the Fourier development of an invariant random field $T$ with respect to the orthonormal basis
$(u_1,\dots, u_d)$
$$
T=\sum_{k=1}^d a_k u_k
$$
then, as, for $t\in\bT$,
$$
T(t^{-1}x)=\sum_{k=1}^d a_k u_k(t^{-1}x)=\sum_{k=1}^d a_k \chi_k(t)u_k(x)
$$
and $T(t^{-1}x)$ and $T(x)$ have the same distribution for every $t\in\bT$, necessarily, for every $k$ such that the action of $\bT$ on $U_k$ is not trivial (that is $\chi_{k}(t)\not\equiv 1$), the coefficient $a_k$ must be invariant in distribution with respect to rotations of the complex plane (and therefore, if it is Gaussian, it must be complex Gaussian).

Remark also that the action of $\bT$ over $V$ cannot be trivial, that is $\chi_k\not\equiv 1$ for some $k$ necessarily.
Actually, as all maximal tori are conjugated (that is if $\bT'$ is another maximal torus then $\bT'=g\bT g^{-1}$ for some $g\in G$) then the action of all maximal tori on $V$ would be trivial which is impossible as the union of all maximal tori is the group itself so that this would imply that the action of $G$ itself is trivial, whereas we assumed $V$ to be irreducible and with dimension $d>1$.

The property, mentioned above, of the random coefficients with respect to the basis of the spherical harmonics in the case of the $\bS^2$, appears now as a particular case.
\end{remark}
\section{Invariant random fields with independent Fourier coefficients}\label{par-complex}
In this section we see results that state that independence assumptions on the Fourier coefficients
implies Gaussianity of the coefficients and of the corresponding random field.

Theorems \ref{caract-complex} and \ref{zero} below are already known (see \cite{MR2342708}) and we reproduce them only to be self-contained, our main concern being the investigation of the validity of Assumption \ref{assum0}, which is a necessary condition in many situations of interest.
\begin{theorem}\label{caract-complex} Assume $G$ to be connected and let $T$ a a.s. square integrable
$G$-invariant random field on the homogeneous space $\cl X$. Let $V$ an irreducible $G$-module
of $L^2(\cl X)$ with dimension $d>1$ and let us assume that coefficients $(a_k)_k$ of the development
$$
T_V=\sum_{k=1}^da_kv_k
$$
with respect to an orthonormal basis $(v_k)_k$ of $V$ are independent. Then they are necessarily Gaussian and the random field $T_{V}$ is Gaussian itself.
\end{theorem}
Remark that in the statement of Theorem \ref{caract-complex}, as in  Theorem \ref{zero} below, we make no assumption concerning the integrability or the existence of finite moments of the r.v.'s $T(x)$ and/or $a_k$. But, of course, under the assumptions of the theorem it follows that necessarily the r.v.'s $T_V(x)$ and $a_k$ has finite moments of every order.

The proof of Theorem \ref{caract-complex} relies on the following Skitovich-Darmois theorem, actually proved in this version by S.~G.~Ghurye and I.~Olkin
\cite{MR0137201} (see also \cite{MR0346969}).
\begin{theorem} \label{Ghurye}Let $X_1,\dots ,X_r$ be mutually
independent random vectors with values in $\R^{n}.$ If, for some real nonsingular $n\times n$ matrices $A_{j},B_{j}, $ $j=1,\dots ,r,$ there are two linear statistics
\begin{equation*}
L_{1}=\sum_{j=1}^{r}A_{j}X_{j}, \qquad L_{2}=\sum_{j=1}^{r}B_{j}X_{j}
\end{equation*}%
that are independent, then the vectors $X_{1},\dots ,X_{r}$ are Gaussian.
\end{theorem}
\noindent{\it Proof} of Theorem \ref{caract-complex}. Let us denote again by $D(g)$ the representative matrix of the left action
of $g\in G$ on $V$ with respect to the orthonormal basis $(v_k)_k$ and by $a$ the vector of the coefficients $a_k$. Thanks to \paref{rotated2} we have
$$
a\enspace\mathop{=}^{distr}\enspace a^g=D(g^{-1})a\ .
$$
Let $1\le k_1<k_2\le\dim V$. Then the joint distribution of $a_{k_1}$ and $a_{k_2}$ is the
same as the joint distribution of
$$
\sum_{j=1}^dD_{k_1j}(g^{-1})a_j \quad\hbox{and}\quad\sum_{j=1}^dD_{k_2j}(g^{-1})a_j
$$
which are therefore themselves independent. Thus we have found two linear statistics of the
r.v.'s $a_k$ that are independent. By the Skitovich-Darmois Theorem \ref{Ghurye} therefore
the joint distribution of the $a_k$'s is Gaussian complex valued as soon as we are able to prove that there exists at least
one element $g\in G$ such that the  real linear transformations
$$
\bC\ni z\to D_{k_1j}(g)z\quad\mbox{ and }\quad \bC\ni z\to D_{k_2j}(g)z,\qquad j=1,\dots,d
$$
are non degenerate. This follows from analyticity properties of the coefficients, as explained at the end of the proof of  Proposition 4.8 of \cite{MR2342708} (here connectedness is required).
\finedim
Remark that the result above does not hold for one-dimensional $G$-modules. As shown in \cite{MR2342708} Example 3.7, it is possible
to construct a non-Gaussian invariant random field on the torus having all its coefficients independent.

Theorem \ref{caract-complex} is not really satisfactory
because its assumptions are not satisfied in the case of
{\it real} random field's, for which the coefficients are necessarily
constrained by the fact that the imaginary parts must cancel and
therefore cannot, in general, be independent. Theorem \ref{caract-complex}
has however its own interest because it contains the essence of the
arguments that we use in the sequel and because it holds without any assumption concerning the
orthonormal basis $(v_k)_k$ of $V$.
\bigskip

\noindent We consider now the case of a {\it real} $G$-invariant random field $T$.

The results are different according to the fact that the irreducible $G$-module $V$ under consideration is of type $\cl I^+$ or $\cl I^0$, as classified at the end of section \ref{par-fourier}.

In the case $V\in\cl I^+$ the fact that we deal with a real random field does not impose constraints on the coefficients of $T_{V}$ and $T_{\overline V}$ but, of course, in order to obtain a real random field one must impose that in the sum $T_{V}+T_{\overline V}$ the imaginary parts cancel.
In this situation the random fields  will be both complex, in general, but their sum will give rise to a real random field.

In this situation Theorem \ref{caract-complex} states that, if $\dim V>1$ and the coefficients $(a_k)_k$ are independent, then $T_{V}$ and $T_{\overline V}$ are both Gaussian and their sum is a Gaussian real random field

Conversely if $V\in\cl I^0$ the coefficients must satisfy some constraints in order to ensure that the random field is real. It is natural then to consider the setting of an  orthonormal basis that is self-conjugated.

It is actually appealing to consider a basis that is formed by real functions. For such a basis, $(v_k)_k$, under the assumption that the coefficients $a_k=\int_{\bS^2}T(x)v_k(x)\, dx$ are independent, Theorem  \ref{caract-complex} applies so that if the random field $T$ is invariant then the $a_k$'s
are jointly Gaussian. Such a statement is not however satisfactory, as we explain below in Remark \ref{real-basis}.

We shall therefore now consider the case of  a basis  $(v_k)_{-\ell\le k\le \ell}$ of $V$ such that
\begin{equation}\label{self-c}
v_{-k}=\overline{v_k}
\end{equation}
This means that 
we assume that, if the dimension of $V$ is odd, the basis contains only one element $v_0$ which is a real function.
If the dimension of $V$ is even we shall still write $(v_k)_{-\ell\le k\le \ell}$ in order to simplify the
notations (there is however no $v_0$ function). In the following arguments we shall consider the case where $\dim V$ is odd, the
case $\dim V$ even being quite similar.

For a basis satisfying \paref{self-c} the fact that $T$ is real imposes to the coefficients the requirement $a_{-k}=\overline{a_k}$. This is the usual setting in the case $\cl X=\bS^2$, where, to be precise, usually one considers the
basis of the spherical harmonics for which it holds $v_{-k}=(-1)^k\overline{ v_k}$ so that the condition above becomes $a_{-k}=(-1)^k\overline{a_k}$, a slight difference that
does not change things. 

The argument in this case can be implemented along the
same lines as in  Theorem \ref{caract-complex}.
Let us assume that the coefficients $(a_k)_{k\ge 0}$, are independent. 
Then, if $m_1\ge 0$, $m_2\ge 0$, $m_1\not=m_2$, and we denote as above by $D_{m,m'}(g)$ the matrix
elements of the action of $G$ on $V$, the two complex r.v.'s
\begin{equation}\label{tilde1}
\displaystyle \widetilde a_{m_{1}}=\sum_{m=-\ell}^\ell D_{m_{1},m}(g^{-1})a_{m}\quad\mbox{and}\quad
\displaystyle \widetilde a_{m_{2}}=\sum_{m=-\ell}^\ell D_{m_{2},m}(g^{-1})a_{m}
\end{equation}
have the same joint distribution as $a_{m_{1}}$ and $a_{ m_{2}}$ and are
therefore independent. The Skitovich-Darmois theorem cannot be applied as before, as the r.v.'s ${a_{m}}$ and
$a_{-m}=\overline{a_{m}}$ are certainly not independent. But \paref {tilde1} can be written
\begin{align*}
\widetilde a_{m_{1}}&=D_{ m_{1},0}(g^{-1})a_{0}+ \sum_{m=1}^\ell
\Bigl(D_{m_{1},m}(g^{-1})a_{m}+D_{ m_{1},-m}(g^{-1})
\overline{a_{m}}\Bigr)\\
\widetilde a_{m_{2}}&=D_{m_{2},0}(g^{-1})a_{0}+ \sum_{m=1}^\ell\Bigl(
D_{m_{2},m}(g^{-1})a_{m}+D_{ m_{2},-m}(g^{-1}) \overline{a_{m}}\Bigr)
\end{align*}
so that $\widetilde a_{m_{1}},\widetilde a_{m_{2}}$ are (real) linear functions of the independent r.v.'s $a_0,\dots, a_\ell$.
Therefore we can again apply Theorem  \ref{Ghurye} as soon as we prove that
$g\in G$ can be chosen so that the real linear applications
\begin{equation}\label{condition}
z\to D_{m_{i},m}(g^{-1})z+D_{m_i,-m}(g^{-1})\overline z, \qquad
m=1,\dots,\ell,\ i=1,2
\end{equation}
are non singular.
It is immediate that this is equivalent to
\begin{equation}\label{diff-moduli}
|D_{m_{i},m}(g^{-1})|\not=|D_{m_i,-m}(g^{-1})|, \qquad m=1,\dots,\ell,\ i=1,2
\end{equation}
We are therefore led to state our main result under the assumption that (\ref{diff-moduli}) is fulfilled.
\begin{assum} \label{assum0} (The mixing condition) Let $V\subset L^2(\cl X)$ a self-conjugated irreducible $G$-module and let
$(v_i)_{-\ell\le i\le \ell}$ a self-conjugated orthonormal basis of $V$ let us denote by $D(g)$ the representative
matrix of the action of $G$ on $V$. We say that $(v_i)_{-\ell\le i\le \ell}$ is
\emph {mixing} if there exist $g\in G$ and $0\le m_1<m_2\le \ell$ such that
\begin{equation}\label{eq-assum0}
|D_{m_{i},m}(g)|\not=|D_{m_{i},-m}(g)|
\end{equation}
for every $0< m\le \ell$, $i=1,2$.
\end{assum}
We have therefore proved the following.
\begin{theorem}\label{zero} Assume $G$ to be connected.
Let $V\subset L^2(\mathcal X)$ be an  irreducible $G$-module such that
$\overline V=V$. Let $(v_k)_{-\ell\le k\le \ell}$ be a self-conjugated mixing (see Assumption \ref{assum0})
basis of $V$. Consider the \emph{real} random field
\begin{equation}\label{rf-zero}
T_{V}(x)=\sum_k a_kv_k(x)
\end{equation}
where the r.v.'s $a_k,k\ge 0$ are independen. Then if $T_{V}$ is $G$-invariant the r.v.'s $(a_k)_{k}$ are jointly Gaussian and
therefore also $T_{V}$ is
Gaussian.
\end{theorem}
\begin{remark}\label{real-imaginary}\rm
It is relevant to point out that in Theorems \ref{caract-complex} and \ref{zero} we do not assume independence of
the real and imaginary parts of the coefficients. Actually under this additional assumption the statement would
become almost trivial (and much weaker) in many situations, as often the invariance of the random field implies
that the coefficients (some of them at least, see Remark \ref{max-torus}) have a distribution that is invariant
with respect to rotations of the complex plane. And it is well known that this assumption together with
independence of the components implies a joint Gaussian distribution, with no need of Assumption \ref{assum0}
(immediate consequence of the Bernstein-Kac theorem as recalled in Proposition \ref{bernstein-kac} below).

This point is important with respect to one of the practical consequences of these results, which is the simulation
of invariant random fields. Actually a natural and computationally efficient procedure in order
to simulate a random field on $\cl X$ is by sampling its Fourier coefficients.
For the case $\cl X=\bS^2$, for instance, Theorem \ref{zero} together with the fact that the basis of the spherical
harmonics is mixing (\cite{dogiocam} p.145) entails that if the coefficients
$a_{\ell m}$'s, $m\ge 0$, of the corresponding development are independent, then, in order to obtain an invariant random field,
they must be Gaussian and the resulting random field will be Gaussian itself. Different choices of their distribution  will lead to a random field which cannot be invariant.
In particular the choice of independent r.v.'s $a_{\ell m}$'s, $m\ge 1$ with a complex Cauchy distribution, for example, cannot produce an invariant random field, {\it even if the real and imaginary parts of $a_{\ell m}$ are not independent}.

Theorem \ref{teo-esseotre} below will imply that if the coefficients
$a_{\ell m}$'s, $m\ge 0$, of the development are independent then, if the resulting random field is invariant, the
$a_{\ell m}$'s are Gaussian, even if the development is made with respect to a generic self-conjugated basis (not necessarily
the spherical harmonics).
\end{remark}
\begin{remark}\label{real-basis}\rm
If the $G$-module $V$ is self-conjugated as in Theorem \ref{zero}, it is natural to consider
an orthonormal basis on $V$ that is formed by real functions.
For instance in the case $\cl X=\bS^2$ one might choose the orthonormal basis given by $v_{\ell 0}=Y_{\ell 0}$ and
$$
v_{\ell m}=\frac 1{\sqrt{2}}\,(Y_{\ell m}+(-1)^mY_{\ell,-m}),\quad v_{\ell, -m}=\frac 1{i{\sqrt{2}}}\,(Y_{\ell m}-(-1)^mY_{\ell,-m}),\quad m\ge 0\ .
$$
The functions $v_{\ell m}$ are real and, if we denote $a_{\ell m}$ the coefficients of the real random field $T$
with respect to the basis of the spherical harmonics, then the coefficients with respect to the basis
$(v_{\ell m})_m$ would be $b_{\ell 0}=a_{\ell 0}$ and
$$
b_{\ell m}=\sqrt{2}\func{Re} a_{\ell m},\enspace  b_{\ell,-m}=\sqrt{2}\func{Im}a_{\ell,m},  \hbox to 0pt{$ \qquad m\ge 1$,}
$$
They are of course real r.v.'s.
A repetition of the arguments of Theorem \ref{caract-complex} now gives immediately that
invariance of the random field and independence of the coefficients $(b_{\ell m})_m$ imply joint Gaussianity of
the coefficients $(b_{\ell m})_m$ without bothering with Assumption \ref{assum0}. This would however be a
much weaker result, as independence of the $(b_{\ell m})_m$'s would imply independence of the real and
imaginary parts of the $(a_{\ell m})_m$'s, which is not required in Theorem \ref{zero} as pointed out above in
Remark \ref{real-imaginary}.
\end{remark}
Let $V$ an irreducible self-conjugated $G$-module, $T_V$ a real random field as in \paref{rf-zero} and
$(v_k)_{-\ell\le k\le \ell}$ a self-conjugated basis as above. Then by Theorem \ref{zero}, under
Assumption \ref{assum0}, if the coefficients $a_k$, $k\ge 0$, with respect to the given basis are
independent they are Gaussian. Moreover, by Corollary \ref{explicative}, as they must have the
same variance and be orthogonal, there exists $c\ge 0$ such that for $k\not=0$ $\E[(\Re a_k)^2]=\E[(\Im a_k)^2]=\frac c2$
(this is a consequence of the orthogonality of $a_k$ and $a_{-k}=\overline{a_k}$) and
 $\E[a_0^2]=c$ (if the basis contains a real function $v_0$).
Conversely, is a real random field $T_V$ with these properties invariant? This
is the object of the next statement.
\begin{theorem}\label{sufficient0} Let $V\subset L^2(\cl X)$ a self-conjugated irreducible $G$-module
and $(v_k)_{-\ell\le k\le \ell}$ a self-conjugated orthonormal basis of $V$
(possibly $k\not=0$ if $\dim H$ is even). Let $T$ a real a.s. square integrable random field on $\cl X$ and
let $(a_k)_{-\ell\le k\le \ell}$ be its random coefficients with respect to the basis above. Then if
the real and imaginary parts of the r.v.'s $a_k, k\ge 0$ (resp. $k>0$ if $\dim V$ is even)
are centered, independent and Gaussian and, for $k\not=0$, there exists $c\ge 0$ such that
$\E[(\Re a_k)^2]=\E[(\Im a_k)^2]=\frac c2$ and $\E[a_0^2]=c$, then the random field
$$
T_V=\sum_{k=-\ell}^\ell a_kv_k
$$
is invariant.
\end{theorem}
\proof We make the proof under the assumption that $\dim V$ is odd, the case of an
even dimension being quite similar. Therefore in the basis $(v_k)_{-\ell\le k\le \ell}$ we have $v_{-k}=\overline {v_k}$
and $v_0$ is a real function. Let $A$ be the matrix of the transformation $\bC^{2\ell+1}\to \bC^{2\ell+1}$
\begin{equation}\label{matrixA}
\begin{pmatrix}
\vphantom{\frac 1{\sqrt{2}}}z_\ell\cr
\vdots\cr
\vphantom{\frac 1{\sqrt{2}}}z_1\cr
z_0
\cr
\vphantom{\frac 1{\sqrt{2}}}z_{-1}\cr
\vdots\cr
\vphantom{\frac 1{\sqrt{2}}}z_{-\ell}
\end{pmatrix}
\qquad\to\qquad
\begin{pmatrix}
\frac 1{\sqrt{2}}\,(z_\ell+z_{-\ell})\cr
\vdots\cr
\frac 1{\sqrt{2}}\,(z_1+z_{-1})\cr
z_0\cr
\frac 1{i\sqrt{2}}\,(z_1-z_{-1})\cr
\vdots\cr
\frac 1{i\sqrt{2}}\,(z_\ell-z_{-\ell})\cr
\end{pmatrix}
\end{equation}
Lemma \ref{desentrelacement} below proves that the matrices $\widetilde D(g)=AD(g)A^{-1}$ are real
orthogonal. Let $a_k=X_k+iY_k$, $k>0$, and $a_0=Z$.
The real r.v.'s $Z,X_k,Y_k$, $k=1,\dots,\ell$, are independent and the matrix $A$
maps the vector $a=(a_\ell,\dots, a_{-\ell})^t$ into
$$
\widetilde a=\begin{pmatrix}
\sqrt{2}\,X_\ell\cr
\vdots\cr
\sqrt{2}\,X_1\cr
Z\cr
\sqrt{2}\,Y_1\cr
\vdots\cr
\sqrt{2}\,Y_\ell\cr
\end{pmatrix}
$$
As the distribution of $\widetilde a$ is Gaussian with all its coordinates centered and independent with
a common variance,
$\widetilde a$ is invariant in distribution under the action of every orthogonal matrix and therefore under the action of
$\widetilde D(g)$ for every $g\in G$. Therefore the random vector $a$ is invariant in distribution under the action of $D(g)$ for every $g\in G$, which
implies the invariance of $T_V$.
\finedim
\begin{lemma}\label{desentrelacement}
$\widetilde D(g)=AD(g)A^{-1}$ is a real orthogonal matrix for every $g\in G$.
\end{lemma}
\proof It is immediate that the rows of $A$ are pairwise orthogonal and unitary. Therefore $A$ is a unitary matrix
as well as $\widetilde D(g)$. Let us prove that $\widetilde D(g)$ maps $\R^{2\ell+1}$ into $\R^{2\ell+1}$, which will end
the proof. Let $\Xi=(\xi_\ell,\dots,\xi_1,\zeta,\eta_1,\dots,\eta_{\ell})^t\in\R^{2\ell+1}$ and
$z_k=\xi_k+i\eta_k$, $k\ge 1$, $z_0=\zeta$. Let $z=(z_{\ell},\dots,z_1,z_0,\overline{z_1},\dots,  \overline{z_{\ell}})^t$; therefore $z=A^{-1}\Xi$ is of the form
\begin{equation}\label{des2}
z=\begin{pmatrix}
\frac1{\sqrt{2}}\,(\xi_\ell+i\eta_\ell)\cr
\vdots\cr
\frac1{\sqrt{2}}\,(\xi_1+i\eta_1)\cr
\zeta\cr
\frac1{\sqrt{2}}\,(\xi_1-i\eta_1)\cr
\vdots\cr
\frac1{\sqrt{2}}\,(\xi_\ell-i\eta_\ell)\cr
\end{pmatrix}
\end{equation}
Now
$$
f(x)=\sum_{k>0} z_kv_k+\zeta v_0+\sum_{k>0} z_{-k}v_{-k}
$$
is a real function. Remark that the matrix $A$ changes any vector of the form \paref{des2} into a
vector of $\R^{2\ell+1}$. As $D(g)z$ is the vector of the coefficients of the function $L_gf$, which is still a real function, its coefficients
are again of the form \paref{des2}, so that $AD(g)A^{-1}\Xi\in\R^{2\ell+1}$.
\finedim
\section{On the validity of the main assumption}\label{par-main}%
In this section we investigate the validity of Assumption \ref{assum0}.

Let us remark first that, for a given self-conjugated $G$-module $V$ of $L^2(\cl X)$, this assumption, as far as we know, might be true for some orthonormal bases of $V$
and not for other ones. So far it is known that it is true for the basis formed by the spherical harmonics when 
$\cl X=\mathbb S^2$ (see \cite{dogiocam} p.144, for a proof), if $\dim V>3$. Actually, as explained below, if 
$\dim V\le 3$ Assumption \ref{assum0} cannot be true. However we investigate the implication between independence of
the coefficients and Gaussianity for the $3$-dimensional irreducible $G$-module of $L^2(\bS^2)$ in Theorem \ref{accadi}.  

Remark that, as $D_{mk}(g)=\langle L_gv_k,v_m\rangle$, condition \paref{eq-assum0}  is equivalent to
\begin{equation}\label{eq-assum20}
|\langle L_gv_m,v_{m_i}\rangle|\not=|\langle L_gv_{-m},v_{m_i}\rangle|\mbox{ for some $g\in G$ and every $0< m\le \ell$, $i=1,2$}\ .
\end{equation}
The main result of this section is Proposition \ref{wedge} where we state a condition equivalent to Assumption
\ref{assum0} carrying a more geometric meaning. This will be the key tool in the next section, where we prove
that every self-conjugated orthonormal basis of an irreducible $G$-module of $L^2(\mathbb S^2)$ with $\dim(V)>3$
is mixing. In \S\ref{par-spheres}
we check the
validity of Assumption \ref{assum0} for the sphere $\cl X=\bS^3$ under the action of $G=SO(4)$, at least for a
class of self-conjugated orthonormal bases.

Let us first state some remarks.
\begin{remark}\label{remarks}\rm
{a)} Mixing (Assumption \ref{assum0}), is a property of the orthonormal self-conjugated basis of the irreducible
$G$-module $V\subset L^2(\cl X)$ and, as far as we know, might hold for some orthonormal self-conjugated basis
and not for other ones. However if it holds for a self-conjugated orthonormal basis $(v_k)_{-\ell\le k\le \ell}$,
then it also holds for every other basis $(w_k)_{-\ell\le k\le \ell}$ of the form $w_k=L_{g_0}v_k$ for some
$g_0\in G$. Actually if $\widetilde D(g)$ denotes the matrix of the action of $G$ on $V$ with respect to the
basis $(w_k)_{-\ell\le k\le \ell}$, that is
$$
w_k(g^{-1}x)=\sum_{i=-d}^d \widetilde D_{ik}(g)w_i(x)
$$
then $\widetilde D(g)=D(g_0^{-1}gg_0)$ so that Assumption \ref{assum0} holds also for $(w_k)_{-\ell\le k\le \ell}$.

b) Assumption \ref{assum0} cannot be true if $v_{m_i}$ is a real function of $L^2(\cl X)$. Actually, as
the left regular action commutes with conjugation,
$$
D_{m_{i},-m}(g)=\langle L_gv_{-m},v_{m_i}\rangle=\langle\overline {L_gv_{m}},v_{m_i}\rangle=\overline{\langle L_gv_{m},v_{m_i}\rangle}=\overline{D_{m_{i},m}(g)}\ .
$$
Therefore $|D_{m_{i},-m}(g)|=|D_{m_{i},m}(g)|$ for every $g\in G$.
This implies that Assumption \ref{assum0} cannot be satisfied if $\overline V=V$ and $\dim V=2$ or $\dim V=3$. Actually in the first case there
is only one $m_i\ge0$, whereas, if $\dim V=3$, the values $m_i=0,1$ are possible, but $v_0$ must be a real function so that \paref{eq-assum0} cannot be satisfied for this value of $m_i$.
The case $\dim V=3$ is of interest because it appears in the Peter-Weyl decomposition of $L^2(\bS^2)$.

In the next section we prove however that also for this $G$-module if the random field is invariant and the coefficients $a_0,a_1$ are independent,
then they are necessarily Gaussian with no need of using the Skitovitch-Darmois theorem and therefore in a situation in which Assumption \ref{assum0} is not satisfied.
This raises the question whether one might prove Theorem \ref{zero} using a different characterization of the Gaussian distribution
than the one provided by
Theorem \ref{Ghurye}, possibly avoiding the need of Assumption \ref{assum0} (see \S\ref{par-conclusions} for a more precise discussion on open questions).
\end{remark}
\vskip-6pt
\finedim
Let $H$ be a irreducible unitary $G$-module of real type (recall Remark \ref{remarks2}) and let $J:H\to H$ a conjugation such that $J^2=1$ and let us denote $\langle\enspace,\enspace\rangle$ the corresponding $G$-invariant scalar product on $H$.
Let $(h_k)_{-\ell\le k\le \ell}$ be a orthonormal basis of $H$ self-conjugated with respect to $J$, i.e. such that $h_{-k}=J(h_k)$. We shall say that such a basis is $J$-mixing if \paref{eq-assum0} holds, now denoting by $D(g)$ the matrix of the action of $G$ on $H$ with respect to this basis. Of course $J$-mixing coincides with mixing if $H\subset L^2(\cl X)$ with the left regular action and $J$ being the usual conjugation $Jv=\overline v$.
\begin{lemma}\label{conjug}
Let $H$ be a irreducible unitary $G$-module of real (resp. quaternionic) type (recall Remark \ref{remarks2}) and let $J:H\to H$ a conjugation such that $J^2=1$ (resp. $J^2=-1$), then, for every $v,w\in H$,
$$
\langle Jv,Jw\rangle=\overline{\langle Jv,Jw\rangle}\ .
$$
\end{lemma}
\proof
It is immediate that
$$
\langle v,w\rangle'=\overline{\langle Jv,Jw\rangle}
$$
is also a $G$-invariant scalar product on $H$, hence, by Schur lemma, there exists a real number $\lambda>0$ such that $\langle v,w\rangle'=\lambda \langle v,w\rangle$. The relation $J^2=1$ (resp. $J^2=-1$) easily implies $\lambda=1$.
\finedim
We are going to express Assumption \ref{assum0} in terms of the action of $G$ on the wedge product $\bigwedge^2 H$.
Recall that
$\bigwedge^2 H$ is endowed with the usual $G$-invariant scalar product
$$
\langle v_1\wedge w_1,v_2\wedge w_2\rangle_2:=\langle v_1,v_2\rangle \langle w_1,w_2\rangle- \langle v_1,w_2 \rangle \langle w_1,v_2\rangle\ .
$$
and we denote $g(v\wedge w)=gu\wedge g w $ the action of $G$ on $\bigwedge^2 H$.
\begin{proposition}\label{wedge}
Let $g\in G$, then $|D_{m_{i},m}(g)|=|D_{m_{i},-m}(g)|$ if and only if
$$
\langle g(h_{m_{i}}\wedge h_{-m_{i}}),h_m\wedge h_{-m}\rangle_2=0\ .
$$
\end{proposition}
\proof We have
$$
\dlines{
\langle g(h_s\wedge h_{-s}),h_m\wedge h_{-m}\rangle_2=
\langle gh_s\wedge gh_{-s},h_m\wedge h_{-m}\rangle_2=\cr
=\langle gh_s,h_m\rangle \langle gh_{-s},h_{-m}\rangle- \langle gh_{s},h_{-m} \rangle \langle gh_{-s},h_{m}\rangle=\cr
=\langle gh_s,h_m\rangle \langle gJh_{s},Jh_{m}\rangle- \langle gh_{s},h_{-m} \rangle \langle gJh_{s},Jh_{-m}\rangle=\cr=|\langle gh_s,h_m\rangle|^2-|\langle gh_{-s},h_{m}\rangle|^2=|D_{m,s}(g)|^2-|D_{m,-s}(g)|^2\cr
}
$$
where we used the fact that $\langle Jv,Jw\rangle=\overline{\langle Jv,Jw\rangle}$ thanks to Lemma \ref{conjug}.
\finedim
Assumption \ref{assum0} can therefore be rephrased in terms of orthogonality of the $G$-orbits of the vectors
 $h_m\wedge h_{-m}$ in $\bigwedge^2 H$.

To be precise, let, for every $1\le m\le\ell$, $W_m\subset \bigwedge^2 H$ the subspace generated by the $G$-orbit of
$h_m\wedge h_{-m}$; let $S$ the set of the pairs $(i,j)$, $1\le i,j\le \ell$, such $W_i\subset W_j^\perp$. Let $\widetilde S$
the set of the indices $1\le i\le \ell$ such that $(i,j)\in S$ for some $1\le j\le \ell$.
\begin{cor}
Assumption \ref{assum0} holds if and only if the complement set ${\widetilde S}^c$ contains at least two indices. In particular
Assumption \ref{assum0} is verified if $\ell\ge 2$ and $S$ is empty.
\end{cor}
\proof
Let $(i,j)$, $1\le i,j\le \ell$ and let
$$
F_{i,j}=\{g\in G; |D_{i,j}(g)|\not=|D_{i,-j}(g)|\}
$$
If $F_{i,j}\not=\emptyset$ then it is a dense open set of $G$. Assumption \ref{assum0} holds if and only if for every
$1\le m_1<m_2\le \ell$ we have that $F_{m_1,m}\not=\emptyset$ for every
$1\le m\le\ell$ and $F_{m_2,m'}\not=\emptyset$ for every
$1\le m'\le\ell$. Now it is sufficient to observe that, by Proposition \ref{wedge},
$F_{i,j}=\emptyset$ if and only if $(i,j)\in S$.
\finedim
In the remainder of this section we introduce a family of orthonormal bases of a $G$-module that arises naturally (the spherical harmonics are of this type) and for which the investigation of the validity of Assumption \ref{assum0} might be simpler.

Let $H$ be an irreducible $G$-module, $\bT$ a maximal torus of $G$ and let us go back to the setting of Remark \ref{max-torus} and consider the decomposition \paref{max-decomposition}. It is possible to assemble
an orthonormal basis of $H$ by picking a unitary vector $u_k$ in each of the $U_k$'s. We say that such a basis is {\it associated to the torus $\bT$}.

If among the $U_k$'s there is only one subspace at most that is associated to a given character of $\bT$,
then the decomposition \paref{max-decomposition} is unique and an associated orthonormal basis
is also unique, up to multiplication of its elements by unitary complex numbers. In this case (i.e. if among the $U_k$'s there is only one subspace at most that is associated to a given
character of $\bT$) we say that $H$ is {\it $\bT$-simple}.

We shall see in the next sections that all irreducible
sub-$G$-modules of $L^2(\bS^2)$ and $L^2(\bS^3)$ are $\bT$ simple with respect to the maximal tori
of $G=SO(3)$ or $G=SO(4)$ respectively.

Let us now suppose that $H$ is a real $G$-module and denote by $J$ a real conjugation. If $u\in U_k$ and $t\in\bT$ we have, denoting
$u\to gu$ the action of $G$,
$$
tu=\chi_k(t)u
$$
for some character $\chi_k$ of $\bT$, so that
\begin{equation}\label{torus-conj}
t\, Ju=J\,tu=\overline{\chi_k(t)}\,Ju=\chi_{-k}(t)\,Ju\ .
\end{equation}
Therefore it is easy to see that an orthonormal basis of $H$ associated to $\bT$ can be chosen in such a way that
it is self-conjugated with respect to $J$. We shall denote by
$(h_k)_{-\ell\le k\le \ell}$ such an orthonormal basis associated to $\bT$, where the index $k$ ranges among 
the corresponding characters of $\bT$ appearing in the decomposition \paref{max-decomposition}. Then it is clear 
that if the relation
\begin{equation}\label{pseudo-assum}
|\langle gh_k,h_{m_i}\rangle|\not=|\langle gh_{-k},h_{m_i}\rangle|,\quad\mbox{for some }g\in G \mbox{ and for every }k\not=0
\end{equation}
holds for the basis $(h_k)_{-\ell\le k\le \ell}$, then it holds also for every other basis that is associated to $\bT$, as two such bases only differ by
multiplication by a unitary complex number.

It is clear that if, in particular, $H\subset L^2(\cl X)$ with the usual conjugation $J$ and \paref{pseudo-assum} holds,
then also Assumption \ref{assum0} holds for the basis $(h_k)_{-\ell\le k\le \ell}$.

It is immediate that if  $H$ is $\bT$-simple then it is also $\widetilde \bT$-simple for any other
maximal torus $\widetilde\bT$. Actually, $\bT$ and $\widetilde\bT$ being conjugated, if
$\widetilde\bT=g^{-1}\bT g$, a basis $(h_k)_{-\ell\le k\le \ell}$ is associated to $\bT$ if and only if $(gh_k)_{-\ell\le k\le \ell}$ is
associated to $\widetilde\bT$. Thanks to Remark \ref{remarks} a), if \paref{pseudo-assum} is
satisfied for a basis associated to $\bT$, then it is also satisfied by all bases associated to $\widetilde\bT$.

The following result states that if an irreducible $G$-module $H$ is $\bT$-simple and satisfies \paref{pseudo-assum}, then the same is true
for every irreducible $G$-module that is equivalent to $H$.
\begin{proposition}\label{equi-tori}
Let $V\subset L^2(\cl X)$ an irreducible $G$-module with $\overline V=V$ and $H$ a $\bT$-simple $G$-module equivalent to $V$. Then also $V$ is $\bT$-simple.
Moreover if \paref{pseudo-assum} is satisfied by the orthonormal bases of $H$ associated to $\bT$, then the same is true for $V$ and every self-conjugated basis associated to a maximal torus of $V$ is mixing.
\end{proposition}
The proof is straightforward.
\section{The sphere $\bS^2$ and related examples}\label{par-sphere-related}
In this section we prove first that, for every irreducible $G$-module, $G=SO(3)$, of dimension $>3$ of
$L^2(\bS^2)$, every self-conjugated orthonormal basis is mixing.
This extends previous results: see \cite{dogiocam} p.144 where this is proved for the basis of the spherical
harmonics. We also give a proof of the fact that the
statement of Theorem \ref{zero} is true for every self-conjugated orthonormal basis
of the irreducible $d$-dimensional $SO(d)$-module of $L^2(\bS^{d-1})$. This covers in particular the case of the three-dimensional $SO(3)$-module of $L^2(\bS^{2})$ a situation
in which we know that Assumption \ref{assum0} is not satisfied (Remark \ref{remarks} b)).

Let us recall that in the Peter-Weyl decomposition of $L^2(\bS^d)$, $d\ge 3$, all the irreducible modules
for the action of $SO(d+1)$ are self-conjugated (see \cite{faraut} pp.196--197), so that, when dealing with a real random field, in order to apply Theorem \ref{zero} the validity of Assumption \ref{assum0} must be checked.

It is well-known that $SO(3)=SU(2)/\{id,-id\}$ so that the irreducible representations of $SO(3)$
are the representations of $SU(2)$ which are trivial on $\{id,-id\}$
(see again \cite{B-D} or \cite{faraut}). The group $G=SU(2)$ acts on the modules $H_\ell$ formed by the homogeneous polynomials in
$2$ complex variables $z_1,z_2$ of degree $\ell$ in the following way: if
$p\in H_\ell$, then, if $z=(z_1,z_2)$,
\begin{equation}\label{repr-su2}
gp(z_1,z_2)=p(az_1-\overline b z_2, bz_1+\overline a z_2)=p(zg)
\end{equation}
where
\begin{equation}\label{descr-su2}
g=\begin{pmatrix}
a&b\cr-\overline b&\overline a
\end{pmatrix}, \qquad a,b\in\C,|a|^2+|b|^2=1
\end{equation}
denotes a generic element of $G=SU(2)$.
The $SU(2)$-modules $H_\ell$ are irreducible and every irreducible $SU(2)$-module is equivalent to $H_\ell$ for some $\ell=0,1,\dots$
The action of $-id$ in these representations is trivial if and only if $\ell$ is even, so that every
irreducible representation of $SO(3)$ is equivalent to $H_\ell$ for some $\ell$ even.
\begin{lemma}\label{lemma-antisym}
Let $P,Q$ be homogeneous polynomials of degree $\ell\ge 1$ in the two complex variables $z_1,z_2$. Let
\begin{equation}\label{def-D(P,Q)}
D(P,Q):=\det\begin{pmatrix}
\displaystyle \frac {\pa P}{\pa z_1}&\displaystyle \frac {\pa P}{\pa z_2}\cr
\displaystyle \frac {\pa Q}{\pa z_1}&\displaystyle \frac {\pa Q^{\vphantom{\sum}}}{\pa z_2}
\end{pmatrix}\ .
\end{equation}
Then $D(P,Q)$ is a polynomial of degree $2\ell-2$ which vanishes if and only if $P=0$ or if $Q=\lambda\,P$ for some $\lambda\in\bC$.
\end{lemma}
We give the proof of Lemma \ref{lemma-antisym} after the following main result.
\begin{theorem}\label{teo-esseotre}
Let $V\subset L^2(G)$, $G=SO(3)$, an irreducible self-conjugated $G$-module
of dimension $2m+1$. Then, if $m>1$, every self-conjugated orthonormal basis of $V$ is mixing.
\end{theorem}
\proof The proof relies on the characterization of Proposition \ref{wedge}. Let $\ell=2m$ and $H_\ell$ as above. Then \paref{def-D(P,Q)} defines a map $(P,Q)\to D(P,Q)$ from $H_\ell\otimes H_\ell\to H_{2\ell-2}$ which is obviously bilinear and antisymmetric. It is also equivariant with respect to the action of $SU(2)$ (acting both on $H_\ell\otimes H_\ell$ and $H_{2\ell-2}$). Actually, denoting $\pa_zP=(\frac {\pa P}{\pa z_1},\frac {\pa P}{\pa z_2})$,
$$
\pa_z(gP)(z)=\pa_z P(g^t\cdot)(z)=(\pa_z P)(g^tz)g^t
$$
and one concludes easily, as $\det g^t=1$. $L(P\wedge Q)=D(P,Q)$ therefore defines a linear equivariant map
$\bigwedge^2H_\ell\to H_{2\ell-2}$, such that $L(P\wedge Q)=0$ if and only if $P\wedge Q=0$
(thanks to Lemma \ref{lemma-antisym}). 

Let us prove that every orthonormal basis self-conjugated with respect to some conjugation
$\widetilde J:H_\ell\to H_\ell$ (i.e. such that
$\widetilde J(f_{-r})=f_r$) is $\widetilde J$-mixing. Let
$(f_{-m},\dots,f_m)$ such a $\widetilde J$-self-conjugated orthonormal basis. If it were not $\widetilde J$-mixing, then by Proposition \ref{wedge}
there would exist $r,s>0$ and two
mutually orthogonal invariant subspaces $U_1,U_2$ of $\bigwedge^2H_\ell$ such that $f_r\wedge f_{-r}\in U_1$,
$f_s\wedge f_{-s}\in U_2$ (recall that we assume $m>1$). As $f_r$ and $f_{-r}$ are orthogonal,
$L(f_r\wedge f_{-r})\not=0$, so that $L$ does not vanish on $U_1$ and by Schur lemma $U_1$ must contain a
$G$-submodule equivalent to $H_{2\ell-2}$. By the same argument
also $U_2$ must contain a $G$-submodule equivalent to $H_{2\ell-2}$, which is not possible, as, by the
Clebsch-Gordan decomposition
(see \cite{dogiocam} \S3.5 or \paref{Clebsch-Gordan} below) the representation $H_{2\ell-2}$ appears only once
in $H_\ell\otimes H_\ell$ and, a fortiori, in $\bigwedge^2H_\ell$.

Therefore, by Proposition \ref{wedge}, the basis $(f_{-m},\dots,f_m)$ is $\widetilde J$-mixing. 

Now let $(v_{-m},\dots,v_m)$ an orthonormal self-conjugated (in the sense of ordinary conjugation, noted $J$)
basis of $V$. The actions of $G=SO(3)$ on $V$ and $H_{2m}$ are equivalent and therefore there exists a map $A:V\to
H_{2m}$ that intertwines the two actions, that is such that $AL_gv=gAv$ for every $g\in G$, $v\in V$. Up to
multiplication by a constant we can assume that $A$ preserves the scalar product. If we note $\widetilde J
f=AJA^{-1}f$, $\widetilde J$ defines a conjugation on $H_{2m}$ with respect to which $f_r=Av_r$ is a
self-conjugated orthonormal basis. By the first part of the proof we know that such a basis is $\widetilde J$-mixing.
Therefore there exists $g\in SO(3)$ such that
$$
|\langle L_g v_r,v_s\rangle_{V}|=|\langle g f_r,f_s\rangle_{H_{2m}}|\not=|\langle g
f_r,f_{-s}\rangle_{H_{2m}}|=|\langle L_g v_r,v_{-s}\rangle_{V}|
$$
and $(v_{-m},\dots,v_m)$ is mixing itself.
\finedim
\noindent{\it Proof} of Lemma \ref{lemma-antisym}. Assume $P\not\equiv0$. If $D(P,Q)=0$ and $\partial_z P \not\equiv 0$, then there exists a function $\lambda:\bC^2\to\bC$ such that, for every $z\in \bC$,
\begin{equation}\label{rel1}
\pa_zQ=\lambda(z)\pa_zP\ .
\end{equation}
Recall Euler formula for homogeneous functions of exponent $\ell$:
$$
\frac {\pa P}{\pa z_1}\, z_1+\frac {\pa P}{\pa z_2}\, z_2=\ell P
$$
and similarly for $Q$, so that from \paref{rel1} we get $Q=\lambda P$.
On the open set $\bC^2 \setminus \Gamma$, where $\Gamma$ is the set of zeros of $P$, we have
$\lambda=\frac QP$, so that, on $\bC^2 \setminus \Gamma$,
$$
\pa_zQ=\frac QP\,\pa_zP\ .
$$
But from $\lambda=\frac QP$ we have also, for $j=1,2$,
$$
\frac {\pa\lambda}{\pa z_j}=\frac 1{P^2}\Bigl(P\,\frac {\pa Q}{\pa z_j}-Q\,\frac {\pa P}{\pa z_j}\Bigr)
=\frac 1{P^2}\Bigl(P\,\Bigl(\frac QP\frac {\pa P}{\pa z_j}\Bigr)-Q\,\frac {\pa P}{\pa z_j}\Bigr)=0\ .
$$
As $\lambda$ is analytic on $\bC^2 \setminus \Gamma$, this implies that 
$Q=const\ P$ on a non empty open set of $\bC^2$ and therefore $Q=const\ P$ everywhere.
\finedim
We address now the question of the three-dimensional irreducible $G$-module of $L^2(\bS^2)$ to which the previous result does not apply
as Assumption \ref{assum0} is not satisfied (see Remark \ref{remarks} b)).

Actually we prove a more general statement.
The key argument is the following classical characterization of the normal distribution.
\begin{proposition}\label{bernstein-kac}
Let $X=(X_1,\dots,X_m)$ a $\R^m$-valued r.v. such that
\tin{a)} the distribution of $X$ is invariant with respect to the action of $SO(m)$;
\tin{b)} there exist $i,j,1\le i,j\le m$ such that $X_i$ and $X_j$ are independent.

Then $X$ is Gaussian.
\end{proposition}
\proof
We can assume for simplicity that $X_1,X_2$ are independent. Let
$$
A=\begin{pmatrix}
\frac 1{\sqrt{2}}&-\frac 1{\sqrt{2}}&0&\dots&0\cr\frac 1{\sqrt{2}}&\hfill\frac 1{\sqrt{2}}&0&\dots&0\cr
0&0&1&\dots&0\cr&&\dots&&\cr0&&\dots&0&1\cr
\end{pmatrix}
$$
Then $A\in SO(d)$ and, by assumption, $X$ and $AX$ have the same distribution. In particular $(X_1,X_2)$ and $\bigl(\frac 1{\sqrt{2}}\,(X_1-X_2),\frac 1{\sqrt{2}}\,(X_1+X_2)\bigl)$  have the same distribution, so that $\frac 1{\sqrt{2}}\,(X_1-X_2)$ and $\frac 1{\sqrt{2}}\,(X_1+X_2)$ are independent.
By the classical Bernstein-Kac characterization of Gaussian measures, $X_1$ and $X_2$ are therefore Gaussian (see \cite{Chaumont-Yor} pp. 74 and 85 for a simple proof). In order to prove joint Gaussianity of $X$, just remark that rotational invariance implies that the characteristic function of $X$ is of the form
$$
\phi_X(\th)=\psi(|\th|),\hbox to0pt{$\qquad\qquad \th\in \R^d$}
$$
for some function $\psi:\R\to\R$. But, as $X_1$ is Gaussian, by choosing $\th=(\th_1,0,\dots,0)$ we have
$$
\psi(|\th|)=\phi_X(\th)=\e^{-\frac12\,\sigma^2\th_1^2}=\e^{-\frac12\,\sigma^2|\th|^2}
$$
where $\sigma^2=\Var(X_1)$, which allows to conclude.
\finedim
Recall that in the Peter-Weyl decomposition of $L^2(\bS^{d-1})$ the smallest irreducible $G$-module besides the constants, $V_d$ say, has always dimension $d$ exactly (see \cite{faraut} p. 197 e.g.).
\begin{theorem}\label{accadi}
Let $(v_{\frac d2},\dots,v_{-1},v_1,\dots, v_{-\frac d2})$ for $d$ even (resp. $(v_{\frac 12(d-1)},\dots,v_0,\dots, v_{-\frac 12(d-1)})$ for odd $d$) a self-conjugated orthonormal basis,
of the $SO(d)$-module $V_d\subset L^2(\mathbb S^{d-1})$. Let $T$ a real invariant random field on $L^2(\bS^{d-1})$ such that its coefficients $(a_{\frac d2},\dots, a_{1})$ (resp. $(a_{\frac 12(d-1)},\dots, a_{0})$) with respect to this basis are independent. Then they are Gaussian.
\end{theorem}%
\proof Let us make the proof for $d$ odd, $d=2m+1$. By assumption we can write $a_k=X_k+iY_k$,
$a_{-k}=X_k-iY_k$ for $1\le k\le d$, $a_0=Z$, where the r.v.'s $Z$, $(X_1,Y_1)$,\dots,$(X_m,Y_m)$ are independent.

Let $a=(a_{-m},\dots,a_0,\dots,a_m)^t$ and denote by $D(g)$, $g\in G$, the matrices of the left regular action of $G=SO(d)$ on the $G$-module $V_d$ with respect to the given orthonormal basis. By assumption the random vectors $a$ and $D(g)a$ have the same distribution for every $g\in G$.
Let now $A$ be the matrix $\C^d\to\C^d$ defined as in \paref{matrixA} (with $m$ replacing by $\ell$)   so that
$$
Aa=\begin{pmatrix}\sqrt{2}\, X_m\cr\vdots\cr
\sqrt{2}\, X_1\cr Z\cr
\sqrt{2}\, Y_1\cr\vdots\cr \sqrt{2}\, Y_m\end{pmatrix}:=\widetilde a
$$
Thanks to Lemma \ref{desentrelacement} the matrices $\widetilde D(g)=AD(g)A^{-1}$ are orthogonal.
Moreover $g\to\widetilde D(g)$ is an irreducible representation of $G$ of dimension $d$, the representations $D$ and $\widetilde D$ being equivalent.
As there is only one $d$-dimensional irreducible representation of $SO(d)$ (up to equivalence), $g\to A ^{-1}D(g)A$ is
equivalent to the natural action of $SO(d)$ on $\bC^d$ and therefore $1=\det(D(g))=\det \widetilde D(g)$ so that $\widetilde D(g)\in SO(d)$ and the image of $\widetilde D$ is $SO(d)$ itself, as the map $g\to \widetilde D(g)$ is injective.

Therefore invariance of the distribution of $a$ with respect to the matrices $D(g)$ entails
invariance of the distribution of $\widetilde a$ with respect to $SO(d)$. As the r.v.'s $X_1$ and $Z$, for instance, are independent,
Proposition \ref{bernstein-kac} implies that $\widetilde a$ is Gaussian, and therefore also $a$.
\finedim
The previous theorem ensures that, if $T$ is an invariant random field on $\bS^2$ and $V$ is an irreducible
$G$-module of $L^2(\bS^2)$ of dimension $3$, independence of the coefficients $a_1$ and $a_0$ with respect
to any self-conjugated orthonormal basis of $V$ entails Gaussianity of $T_V$, even if Assumption \ref{assum0} is not true for such $V$.
\begin{remark}\rm
If $V_d$ is the $d$-dimensional irreducible $SO(d)$-module and if $d\ge 4$, then Theorem \ref{accadi}
actually follows from
Theorem \ref{zero} as, as a consequence of Proposition \ref{wedge}, every self-conjugated orthonormal basis
is mixing. In fact for $d\ge 5$
the module $\bigwedge^2V_d$ is irreducible whereas, for $d=4$ $\bigwedge^2V_d$,  has two irreducible components
which are the eigenspaces of the Hodge$^*$ operator. Therefore a non zero real vector of the form $iv\wedge \overline v$ cannot be
contained in either eigenspace (see \cite{B-D} pp.272--274).
\end{remark}
\begin{example}\label{so3-su2}\rm ($SO(3)$ and $SU(2)$)
In the same line of arguments it is easy to check that, for a real invariant random field, independence of the coefficients entails
Gaussianity in the cases $\cl X=G=SO(3)$ and $\cl X=G=SU(2)$.

Actually if $\cl X=G=SO(3)$ this is partially known when considering the basis given by the normalized columns
(or rows) of the Wigner matrices: in every isotypical submodule one of the columns is generated by
the spherical harmonics for which it is known that, for $\ell>1$, Assumption \ref{assum0} holds so that Theorem \ref{zero}
applies. As for the other columns, they
are not self-conjugated but conjugated pairwise, so that one can  apply Theorem \ref{caract-complex}.

However Theorems \ref{teo-esseotre} and \ref{accadi} ensure that, even considering a different decomposition
of the isotypical spaces, it is not possible to simulate an invariant non Gaussian random field using independent coefficients.

This is true also for $\cl X=G=SU(2)$ as in the Peter-Weyl decomposition,
in addition to those already
considered for $G=SO(3)$, other representations appear that are quaternionic, so that the corresponding isotypical
modules cannot contain self-conjugated irreducible modules (recall Remark \ref{remarks2}).
\end{example}
\section{The sphere $\bS^3$}\label{par-spheres}
In this section we prove that for every irreducible $G$-module, $G=SO(4)$, of $L^2(\bS^3)$ every basis adapted to a maximal torus is mixing.

We shall need some known facts about the group $SO(4)$ and its representations.
$G=SO(4)$ is isomorphic to $SU(2)\times SU(2)/\{(id,id),(-id,-id)\}$. Therefore its irreducible representations are of the form
$H_\ell\otimes H_k$, $H_\ell,H_k$ being the irreducible modules of $SU(2)$ introduced at the beginning of \S\ref{par-sphere-related},
with the condition
that the action of $(-id,-id)$ is trivial. As these modules are formed by the homogeneous polynomials of degree $\ell$ and $k$ respectively in the complex variables $z_1,z_2$, one has
$$
(-id,-id)(p\otimes q)=(-1)^{\ell+k}p\otimes q
$$
and therefore the irreducible modules of $SO(4)$ are of the form $H_\ell\otimes H_k$ with $\ell+k$ even.
In order to determine
the Peter-Weyl decomposition of $L^2(\cS^3)$, $\cS^3=SO(4)/SO(3)$, one must recall that, in the
isomorphism $G\simeq SU(2)\times SU(2)/\{(id,id),(-id,-id)\}$, $SO(3)$ is mapped into the diagonal.
Then the action of $SO(3)$ on $H_\ell\otimes H_k$ is  $g(p\otimes q)=gp\otimes gq$ of $SU(2)$.
By the Clebsch-Gordan formula for $SU(2)$, the  action of $SU(2)$ on the tensor product $H_\ell\otimes H_k$ can be decomposed as
\begin{equation}\label{Clebsch-Gordan}
H_\ell\otimes H_k=\bigotimes_{j=0}^{d_q}H_{\ell+k-2j},\quad d_q=\min(\ell,k)\ .
\end{equation}
Therefore the trivial representation appears in this decomposition if and only if
$$
\ell+k\mbox{ is even},\quad \frac {\ell+k}2\le \ell,\quad \frac {\ell+k}2\le k
$$
that is if and only if $\ell=k$. We have therefore found that the representations of $SO(4)$
appearing in the Peter-Weyl decomposition of $L^2(\bS^3)$ are exactly those that are equivalent to $H_\ell\otimes H_\ell$.
Remark that the smallest dimension of these, besides the case $\ell=1$ of the constants, is $4$, so that
we do not have to bother with the problem of dimension $3$ appearing for the sphere $\bS^2$, as discussed
in Remark \ref{remarks} b).

On the $SU(2)$-module $H_\ell$ introduced in \S5 let us consider the polynomials $p_s(z_1,z_2)=z_1^sz_2^{\ell-s}$,
$s=0,\dots,\ell$ which form an orthogonal basis with respect to the  scalar product
\begin{equation}\label{inv-scal}
\langle p_s,p_r\rangle =\frac {s!(\ell-s)!}{\ell!}\,\delta_{s,r}=\frac 1{{l\choose s}}\,\delta_{s,r}
\end{equation}
which turns out to be $SU(2)$-invariant. Therefore the polynomials $ e_s=c_s p_s$, $s=0,\dots,\ell$ with $c_s=\sqrt{{\ell \choose s}}$ form an orthonormal basis of the unitary $SU(2)$-module $H_\ell$.

A maximal torus of $SU(2)$ is the subgroup of the elements
$$
t_\th=\begin{pmatrix}\e^{i\th}&0\cr 0&\e^{-i\th}\cr\end{pmatrix}
$$
whose action on the polynomials $p_s(z_1,z_2)=z_1^sz_2^{\ell-s}$ is
$$
t_\th p_s=\e^{i(2s-\ell)\th}\, p_s\ .
$$
Thus, with respect to the invariant scalar product \paref{inv-scal}, the elements $e_s=c_s p_s$ with
$c_s=\sqrt{\ell\choose s}$ form an orthonormal basis of $H_\ell$ that is adapted to the maximal torus $\bT$. In particular
$H_\ell$ is $\bT$-simple.

The following computation is our key argument.
 We have
$$
\dlines{
ge_s=c_s(az_1-\overline b z_2)^{s}(bz_1+\overline a z_2)^{\ell-s}
=c_s\sum_{r=0}^{\ell}\underbrace{\sum_{{{\scriptstyle h+k=r\atop \scriptstyle 0\le h\le s}\atop \scriptstyle 0\le k\le \ell-s}}
{s\choose h}{\ell-s\choose k}a^h\overline
a^{\ell-s-k}b^k(-\overline b)^{s-h}}_{:=H_{r,s}}z_1^rz_2^{2m-r}=\cr
=c_s\sum_{r=0}^{\ell}H_{r,s}\,z_1^{r}z_2^{\ell-r}=c_s\sum_{\ell=-m}^{m}\frac {1}{c_r}\,H_{r,s}\,e_r\cr
}
$$
and therefore
$$
\langle ge_s,e_j\rangle=\frac {c_s}{c_j}H_{j,s}=
\frac {c_s}{c_j}\sum_{{{\scriptstyle h+k=j\atop \scriptstyle 0\le h\le s}\atop \scriptstyle 0\le k\le \ell-s}}
{s\choose h}{\ell-s\choose k}a^h\,\overline
{a}^{\ell-s-k}\,b^k\,\overline b^{s-h}(-1)^{s-h}\ .
$$
Taking into account the condition $h+k=j$ this can be written
$$
\dlines{
\langle ge_s,e_j\rangle=\frac {c_s}{c_j}H_{j,s}=
\frac {c_s}{c_j}\sum_{{{
\scriptstyle 0\le h\le s}\atop \scriptstyle 0\le j-h\le \ell-s}}
{s\choose h}{\ell-s\choose j-h}a^h\>\overline
a^{\ell-s-j+h}\>b^{j-h}\>\overline b^{s-h}(-1)^{s-h}=\cr
=\frac {c_s}{c_j}\,\overline a^{\ell-s-j}\>\overline b^{s-j}\sum_{{{
\scriptstyle 0\le h\le s}\atop \scriptstyle 0\le j-h\le \ell-s}}
{s\choose h}{\ell-s\choose j-h}|a|^{2h}|b|^{2(j-h)}(-1)^{s-h}\cr
}
$$
and therefore
\begin{equation}\label{main-poly}
\begin{array}{c}
\displaystyle|\langle ge_s,e_j\rangle|^2=\cr
\cr
\displaystyle=\frac {c_s^2}{c_j^2}|a|^{2(\ell-s-j)}|b|^{2(s-j)}\biggl(\sum_{{{
\scriptstyle 0\le h\le s}\atop \scriptstyle 0\le j-h\le \ell-s}}
{s\choose h}{\ell-s\choose j-h}|a|^{2h}|b|^{2(j-h)}(-1)^{s-h}\biggr)^2=\cr
\cr
=\displaystyle\frac {c_s^2}{c_j^2}\biggl(\sum_{{{
\scriptstyle 0\le h\le s}\atop \scriptstyle 0\le j-h\le \ell-s}}
{s\choose h}{\ell-s\choose j-h}|a|^{\ell-s-j+2h}|b|^{s+j-2h)}(-1)^{s-h}\biggr)^2:=\cl P^\ell_{s,j}(|a|,|b|)
\end{array}
\end{equation}
which is a homogeneous polynomial of degree $2\ell$ in the variables $|a|,|b|$.
Let us point out that in the sum inside the square defining $\cl P^\ell_{s,j}$ the range of $h$ is
\begin{equation}\label{highest}
\max(0,-\ell+s+j)\le h\le \min(s,j)\ .
\end{equation}
Let us assume first $\ell$ even, $\ell=2m$. Let
$$
f_k=e_{m+k}\ .
$$
$(f_k)_{-m\le k\le m}$ is also an orthonormal basis with respect to the $SU(2)$-invariant scalar product
\paref{inv-scal} and adapted to $\bT$.
The maximal torus of $SU(2)\times SU(2)$ is $\bT\times \bT$. As
$$
(t_{\th_1}, t_{\th_2})f_{k_1}\otimes f_{k_2}=\underbrace{\e^{i(2k_1\th_1+2k_2\th_2)}}_{\chi_{k_1,k_2}(t_{\th_1}, t_{\th_2})}\,f_{k_1}\otimes f_{k_2},\qquad -m\le k_1,k_2\le m
$$
$H_\ell\otimes H_\ell$ is simple with respect to $\bT\times\bT$ and the basis $(f_{k_1}\otimes f_{k_2})_{k_1,k_2}$ is adapted
to the maximal torus above. Moreover if $f_{k_1}\otimes f_{k_2}$ is the eigenvector of the character $\chi_{k_1,k_2}$
of $\bT\times \bT$, then $f_{-k_1}\otimes f_{-k_2}$ is an eigenvector of $\overline \chi_{k_1,k_2}$.
We proceed now to check condition \paref{pseudo-assum} in view of taking advantage of Proposition \ref{equi-tori}.
We must show that for some $m_1>0$, $m_2>0$
\begin{equation}\label{different}
|\langle (g_1,g_2)(f_{m_1}\otimes f_{m_2}),f_{r_1}\otimes f_{r_2}\rangle|^2\not=
|\langle (g_1,g_2)(f_{m_1}\otimes f_{m_2}),f_{-r_1}\otimes f_{-r_2}\rangle|^2\mbox{ for every} -m\le r_1,r_2\le m
\end{equation}
for some $g_1,g_2\in SU(2)$. We have
$$
\dlines{
|\langle (g_1,g_2)(f_{m_1}\otimes f_{m_2}),f_{r_1}\otimes f_{r_2}\rangle|^2=
|\langle g_1f_{m_1},f_{r_1}\rangle|^2|\langle g_2f_{m_2},f_{r_2}\rangle|^2\cr
}
$$
and taking into account \paref{main-poly}
$$
|\langle gf_{k},f_{r}\rangle|^2=|\langle ge_{k+m},e_{r+m}\rangle|^2=P^\ell_{k+m,r+m}(|a|,|b|)
$$
so that, denoting by $a_1,b_1$ and $a_2,b_2$ the coordinates of $g_1$ and $g_2$ in the representation \paref{repr-su2},
\begin{equation}\label{piu}
|\langle (g_1,g_2)(f_{m_1}\otimes f_{m_2}),f_{r_1}\otimes f_{r_2}\rangle|^2
=\cl P^{2m}_{m+m_1,m+r_1}(|a_1|,|b_1|)
\cl P^{2m}_{m+m_2,m+r_2}(|a_2|,|b_2|)
\end{equation}
and
\begin{equation}\label{meno}
|\langle (g_1,g_2)(f_{m_1}\otimes f_{m_2}),f_{-r_1}\otimes f_{-r_2}\rangle|^2
=\cl P^{2m}_{m+m_1,m-r_1}(|a_1|,|b_1|)
\cl P^{2m}_{m+m_2,m-r_2}(|a_2|,|b_2|)\ .
\end{equation}
In order to conclude  we must prove that for some values of $a_1,b_1,a_2, b_2$
with $|a_1|^2+|b_1|^2=|a_2|^2+|b_2|^2=1$ the right-hand sides in \paref{piu} and \paref{meno} are different.
For every $r\not=0$ if the two polynomials $\cl P^\ell_{m+m_i,m-r}$ and $\cl P^\ell_{m+m_i,m+r}$, both
homogeneous of degree $4m$, coincide on the circle $|a|^2+|b|^2=1$, they would coincide on the
whole of $\R^2$. In order to see that this cannot happen we look at the monomial that exhibits the highest
exponent in $|a|$ and see that the degrees are different. Recalling \paref{highest}, we must show that the two values
$$
h_1=\min(m+m_i,m+r)\quad \mbox{and}\quad h_2=\min(m+m_i,m-r)
$$
are different. This is done by checking directly all possibilities: as $m_i>0$, then
$$
\begin{tabular}{|c|c|c|}
\hline
 & $h_1$ & $h_2$ \\
\hline
$0<r\le m_i$ & $m+m_i$ & $m-r$ \\
\hline
$m_i<r\le m$ & $m+r$ & $m-r$ \\
\hline
$-m_i\le r<0$ & $m+r$ & $m+m_i$ \\
\hline
$-m\le r<-m_i$ & $m-r$ & $m+r$ \\
\hline
\end{tabular}
$$
Therefore, unless $r=0$ of course, $h_1\not=h_2$ in all possible occurrences. Therefore
the two polynomials at the right-hand side of \paref{piu} and \paref{meno} are different
if one at least between $r_1$ and $r_2$ is different from $0$.

Along the same lines goes the proof for $\ell$ odd.
Thanks to Proposition \ref{equi-tori}, we have
\begin{theorem}\label{s3} Let $V_\ell$ a
irreducible $G$-module appearing in the Peter-Weyl decomposition of $L^2(\bS^3)$ of dimension $>1$. Then
every self-conjugated basis of $V_\ell$ associated to a maximal torus is mixing.
\end{theorem}
\section{Some open questions}\label{par-conclusions}
This paper gives some precisions about properties of the Fourier coefficients of an invariant random field and clarifies some important points in the direction of characterizing the random fields on an homogeneous space whose coefficients in their Fourier development are independent (or at least that can be simulated through the generation of independent r.v.'s) on the track of \cite{BM06} and \cite{MR2342708}.

However it also points out some natural questions that remain open to conjecture. We make here a tentative list.

\tin{1)} In order to prove the Gaussianity of such random field's we used the Skitovitch-Darmois theorem whose application in turn requires, in many cases of interest, to ascertain that Assumption \ref{assum0} is verified. But we have also remarked that Gaussianity remains true in situations where Assumption \ref{assum0} is not true (see Remark \ref{remarks} b). So one might think of taking advantage of a characterization of Gaussianity different form the one that is provided by the Skitovitch-Darmois theorem, and thus be ridden of Assumption \ref{assum0}.

\tin {2)} It is nevertheless of interest also to investigate the validity of Assumption \ref{assum0}. Is it always
true (at least for real $G$-modules of dimension $> 3$)? We do not know of counterexamples so far.

\tin {3)} In a less ambitious perspective, is Assumption \ref{assum0} true for the groups $SO(d)$ and for the
spheres $\bS^{d-1}$, $d\ge 5$? For which classes of orthonormal bases? Intuition should point towards the
positive: as these structures contain $SO(3)$ and $SO(4)$ for which the result is proved, it might be possible to develop a sort of induction procedure on $d$ in order to obtain the validity of Assumption \ref{assum0} or at least that the Gaussian characterization holds. Remark that one possible way of attacking this problem is through an extension of Proposition \ref{bernstein-kac}: does its statement
remain true if the assumption of invariance with respect to the group $SO(n)$ is replaced by invariance with respect to a subgroup of $SO(n)$ that acts irreducibly on $\R^n$?

\bibliography{bibbase}

\providecommand{\bysame}{\leavevmode\hbox to3em{\hrulefill}\thinspace}
\providecommand{\MR}{\relax\ifhmode\unskip\space\fi MR }
\providecommand{\MRhref}[2]{%
  \href{http://www.ams.org/mathscinet-getitem?mr=#1}{#2}
}
\providecommand{\href}[2]{#2}
\begin{thebibliography}{10}

\bibitem{BM06}
P.~Baldi and D.~Marinucci, \emph{Some characterizations of the spherical
  harmonics coefficients for isotropic random rields}, Statist. Probab. Letters
  \textbf{77} (2007), 490--496.

\bibitem{MR2342708}
P.~Baldi, D.~Marinucci, and V.~S. Varadarajan, \emph{On the characterization of
  isotropic {G}aussian fields on homogeneous spaces of compact groups},
  Electron. Comm. Probab. \textbf{12} (2007), 291--302 (electronic).

\bibitem{B-D}
Th. Br{\"o}cker and T.~tom Dieck, \emph{Representations of compact {L}ie
  groups}, Graduate Texts in Mathematics, vol.~98, Springer-Verlag, New York,
  1995.

\bibitem{Chaumont-Yor}
L.~Chaumont and M.~Yor, \emph{Exercises in probability}, Cambridge Series in
  Statistical and Probabilistic Mathematics, vol.~13, Cambridge University
  Press, Cambridge, 2003.

\bibitem{faraut}
J.~Faraut, \emph{Analysis on {L}ie groups}, Cambridge Studies in Advanced
  Mathematics, vol. 110, Cambridge University Press, Cambridge, 2008.

\bibitem{MR0137201}
S.~G. Ghurye and I.~Olkin, \emph{A characterization of the multivariate normal
  distribution}, Ann. Math. Statist. \textbf{33} (1962), 533--541.

\bibitem{MR0346969}
A.~M. Kagan, Yu.~V. Linnik, and C.~Radhakrishna Rao, \emph{Characterization
  problems in mathematical statistics}, John Wiley \& Sons, New
  York-London-Sydney, 1973.

\bibitem{bib:M}
A.~Malyarenko, \emph{Invariant random fields in vector bundles and application
  to cosmology}, Ann. Inst. H. Poincar\'e Probab. Stat. \textbf{47} (2011),
  no.~4, 1068--1095.

\bibitem{dogiocam}
D.~Marinucci and G.~Peccati, \emph{Random {F}ields}, London Mathematical
  Society Lecture Note Series 389, Cambridge University Press, Cambridge, 2011.

\bibitem{mp-continuity}
\bysame, \emph{Mean-square continuity on homogeneous spaces of compact groups},
  preprint (2012).

\bibitem{Pepy}
G.~Peccati and J.-R. Pycke, \emph{Decompositions of stochastic processes based
  on irreducible group representations}, Teor. Veroyatn. Primen. \textbf{54}
  (2009), no.~2, 304--336.

\bibitem{MR1143783}
N.~Ja. Vilenkin and A.~U. Klimyk, \emph{Representation of {L}ie groups and
  special functions. {V}ol. 1}, Mathematics and its Applications (Soviet
  Series), vol.~72, Kluwer Academic Publishers Group, Dordrecht, 1991.

\end{thebibliography}
\bibliographystyle{amsplain}
\end{document}